\documentclass{amsart}
\usepackage[francais]{babel}
\usepackage{amsmath}
\usepackage{amssymb}
\usepackage{amscd}
\usepackage{amsthm}
\usepackage[T1]{fontenc}
\usepackage{multirow}

\newcommand{\CL}{\mathcal{L}}
\newcommand{\CH}{\mathcal{H}}

\newcommand{\cmp}[2]{{}^{\langle{#2}\rangle\!}{#1}}

\newcommand{\dS}{d_{\mathrm{S}}}
\newcommand{\Gh}{\hat G}
\newcommand{\Lh}{\hat L}

\newcommand{\orb}{\mathcal{O}}
\newcommand{\N}{\mathcal{N}}
\newcommand{\F}{\mathbb{F}}
\newcommand{\Qlb}{\bar{\mathbb{Q}}_\ell}
\newcommand{\Nb}{\bar\N}
\newcommand{\bA}{\bar A}
\newcommand{\bs}{\bar s}
\newcommand{\Ab}{\bar{A}}

\newcommand{\E}{\mathcal{E}}
\newcommand{\Fa}{\mathcal{F}}
\newcommand{\G}{\mathcal{G}}
\newcommand{\M}{\mathcal{M}}

\newcommand{\ta}{\tilde{a}}

\newcommand{\bc}{\mathbf{c}}
\newcommand{\bd}{\mathbf{d}}
\DeclareMathOperator{\ad}{ad}
\DeclareMathOperator{\der}{der}
\DeclareMathOperator{\ind}{ind}
\DeclareMathOperator{\Irr}{Irr}
\DeclareMathOperator{\IC}{IC}
\DeclareMathOperator{\Cl}{Cl}
\DeclareMathOperator{\Ind}{Ind}

\DeclareMathOperator{\Tr}{Tr}
\DeclareMathOperator{\Norm}{N}
\DeclareMathOperator{\Cent}{C}
\DeclareMathOperator{\Centre}{Z}
\DeclareMathOperator{\ord}{ord}
\DeclareMathOperator{\pr}{pr}

\DeclareMathOperator{\sgn}{\varepsilon}
\DeclareMathOperator{\supp}{supp}
\DeclareMathOperator{\unif}{unif}
\DeclareMathOperator{\unip}{uni}

\DeclareMathOperator{\VM}{vm}
\DeclareMathOperator{\VMP}{vmp}
\DeclareMathOperator{\Lu}{L}
\DeclareMathOperator{\GM}{GM}
\def\isom{{\overset\sim\to}}
\def\ie{{\it i.e.,\ }}

\newtheorem{thm}{Th\'eor\`eme}[section]
\newtheorem{prop}[thm]{Proposition}
\newtheorem{lem}[thm]{Lemme}
\newtheorem{cor}[thm]{Corollaire}
\theoremstyle{definition}
\newtheorem{defn}[thm]{D\'efinition}
\newtheorem{remark}[thm]{Remarque}
\theoremstyle{remark}

\catcode`\@=11
\def\courriel#1#2{\begingroup
    \@ifnotempty{#2}{\nobreak\indent{\itshape Adresse courriel}%
      \@ifnotempty{#1}{, \ignorespaces#1\unskip} :\space
      \ttfamily#2\par}\endgroup}
\renewcommand{\email}[2][]{\g@addto@macro\addresses{\courriel{#1}{#2}}}
\catcode`\@=12
\newenvironment{symbole}{\begingroup\setlength{\arraycolsep}{0pt}%
  \left(\begin{matrix}}{\end{matrix}\right)\endgroup}
\newenvironment{ssymbole}{\begingroup\setlength{\arraycolsep}{0pt}%
  \left(\begin{smallmatrix}}{\end{smallmatrix}\right)\endgroup}

\title{Supports unipotents de faisceaux caract\`eres}
\author{Pramod N.~Achar}
\thanks{Le premier auteur \'etait partiellement appuy\'e par une bourse
post-doctorale de la NSF}
\address{Department of Mathematics\\
  University of Chicago\\
  Chicago, IL \ 60637, USA}
\email{pramod@math.uchicago.edu}
\author{Anne-Marie Aubert}
\address{Institut de Math\'ematiques de Jussieu\\
Unit\'e Mixte Paris 6 / Paris 7 CNRS de Recherche 7586\\
Universit\'e Pierre et Marie Curie, F-75252 Paris Cedex 05
  \\
  France}
\email{aubert@math.jussieu.fr}
\date{23 avril 2003}
\begin{document}

\begin{abstract}
Soit $G$ un groupe alg\'ebrique r\'eductif sur la cl\^oture
alg\'ebrique d'un corps fini $\F_q$ et d\'efini sur ce dernier. 
L'existence du support unipotent d'un
caract\`ere irr\'eductible du groupe fini $G(\F_q)$, ou d'un faisceau
caract\`ere de $G$, a \'et\'e \'etablie dans diff\'erents cas par
Lusztig, Geck et Malle, et le second auteur. Dans cet article, nous
d\'emontrons que toute classe unipotente sur laquelle la restriction
du faisceau caract\`ere ou du caract\`ere donn\'e 
est non nulle est contenue l'adh\'erence de Zariski de son support unipotent. 
Pour \'etablir ce r\'esultat, nous \'etudions certaines repr\'esentations
des groupes de Weyl, dites ``bien support\'ees''.
\end{abstract}

\maketitle

\section{Introduction}
\label{sect:intro}

Soit $\bar\F_q$ la cl\^oture alg\'ebrique d'un corps fini $\F_q$, et
soit $p$ sa caract\'eristique.  Soit $G$ un groupe alg\'ebrique
r\'eductif sur $\bar\F_q$ qui est d\'efini sur $\F_q$, et soit $F$
l'endomorphisme de Frobenius associ\'e \`a la structure
$\F_q$-rationnelle de $G$. \'Etant donn\'e un caract\`ere
irr\'eductible $\rho$ de $G^F$, Lusztig a pos\'e en 1980 le probl\`eme
de lui associer une classe unipotente $F$-stable $\orb$ de $G$ telle que
la restriction de $\rho$ \`a $\orb^F$ soit non nulle, et telle que la
dimension de $\orb$ soit maximale parmi les classes unipotentes poss\'edant 
cette propri\'et\'e.  Une telle classe est appel\'ee le \emph{support
unipotent} du caract\`ere.  Lusztig a r\'esolu lui-m\^eme ce
probl\`eme en 1992 \cite{LUS}, sous l'hypoth\`ese que $p$ est
suffisamment grand. Dans le m\^eme article, il a aussi \'etabli un
r\'esultat concernant le support unipotent des faisceaux
caract\`eres. Geck et Malle ont r\'eussi \`a enlever l'hypoth\`ese sur
$p$ dans le cas d'une notion l\'eg\`erement diff\'erente de support
unipotent de caract\`ere \cite{GVM}, \cite{GMU}.  Quant aux faisceaux
caract\`eres, le second auteur a \'etendu le r\'esultat de Lusztig au
cas o\`u la caract\'eristique est bonne \cite{Au}, en utilisant la
m\^eme approche que dans \cite{LCV}.

Le but du pr\'esent article est de raffiner la description de
l'ensemble des classes unipotentes sur lesquelle la restriction d'un
caract\`ere ou d'un faisceau caract\`ere est non nulle, par une
\'etude d\'etaill\'ee des repr\'esentations induites des groupes de
Weyl. Si $\orb$ est une classe unipotente de $G$, soit $A(\orb)$ le
groupe des composantes du centralisateur d'un \'el\'ement de $\orb$.
Soit $W$ le groupe de Weyl de $G$, et soit $\N_G$ l'ensemble des
couples $(\orb,\pi)$ o\`u $\orb$ est une classe unipotente et $\pi$
est une repr\'esentation irr\'eductible du groupe $A(\orb)$. Rappelons
que la correspondance de Springer associe \`a toute repr\'esentation
irr\'eductible de $W$ un \'el\'ement de $\N_G$. Cette application,
not\'ee $\nu\colon \Irr(W) \to \N_G$, est injective. \`A la
section~\ref{sect:weyl}, nous introduisons le concept d'\^etre
``(sp\'ecialement) bien support\'e'': une repr\'esentation de $W$ est
dite (sp\'ecialement) bien support\'ee si les \'el\'ements de $\N_G$
correspondant \`a ses composantes irr\'eductibles v\'erifient
certaines conditions.  Le fait qui rend ce concept utile est le
Th\'eor\`eme~\ref{thm:b-s}, qui affirme qu'une repr\'esentation
induite d'une repr\'esentation sp\'ecialement bien support\'ee est
bien support\'ee. La plus grande partie de la preuve de cet \'enonc\'e
est achev\'ee \`a la section~\ref{sect:weyl}, quoique certains calculs
pour les groupes classiques soient report\'es \`a la
section~\ref{sect:biensupp}. 

La section~\ref{sect:ab} traite des liens entre les d\'eveloppements
des deux sections pr\'ec\'edentes et l'ordre partiel introduit par le
premier auteur en \cite{A} sur les classes de conjugaison du quotient
de Lusztig associ\'e \`a une classe unipotente. Aux
sections~\ref{sect:su-fc}--\ref{sect:su-car} sont \'etablis les
r\'esultats principaux de l'article. Ces \'enonc\'es affirment que
toute classe unipotente sur laquelle la restriction d'un faisceau
caract\`ere ou d'un caract\`ere irr\'eductible, respectivement, est
non nulle est contenue dans l'adh\'erence de Zariski du support
unipotent de celui-ci. La version qui traite des caract\`eres \'etend un
r\'esultat \'etabli par Geck et Malle pour les caract\`eres
unipotents.

Nous notons le normalisateur (resp.~le centralisateur) d'un
sous-groupe $H$ dans un groupe $G$ par $\Norm_G(H)$
(resp.~$\Cent_G(H)$). Le centre de $G$ est not\'e $\Centre(G)$. Si $K$
est un groupe fini, $\Irr(K)$ d\'esignera l'ensemble des classes
d'isomorphismes de repr\'esentations irr\'eductibles de $K$, et
$\Cl(K)$ l'ensemble des classes de conjugaison de $K$.

Nous aimerions remercier F.~L\"ubeck et K.~McGerty pour des
conversations utiles.

\section{Rappels sur les classes unipotentes et les groupes de Weyl}
\label{sect:weyl}

Dans cette section, nous allons rappeler certains concepts cl\'es qui
nous seront utiles dans la suite. Le premier paragraphe traite les
degr\'es g\'en\'erique et fant\^ome d'une repr\'esentation de $W$ et
puis introduit le quotient de Lusztig de $A(\orb)$. Au deuxi\`eme
paragraphe nous r\'evisons la relation entre les cellules bilat\`eres
de $W$ et les pi\`eces sp\'eciales de la vari\'et\'e unipotente. Le
dernier paragraphe est consacr\'e \`a l'\'etude de plusieurs sortes
d'induction. On ne pr\'etend pas d'avoir trait\'e ces th\`emes d'une
mani\`ere approfondie: nous nous contentons ici de simplement donner
une br\`eve liste des faits dont nous aurons besoin plus tard.

\subsection{Le quotient de Lusztig}

\`A toute repr\'esentation d'un groupe de Weyl sont associ\'es deux
polyn\^omes d'une variable $q$, le degr\'e g\'en\'erique et le degr\'e
fant\^ome. Le degr\'e g\'en\'erique donne la dimension de la
repr\'esentation unipotente correspondante du groupe de Chevalley fini
$G(\mathbb{F}_q)$, tandis que le degr\'e fant\^ome est tel que le
coefficient de $q^i$ donne la multiplicit\'e de la repr\'esentation
donn\'ee dans la $i$-\`eme puissance sym\'etrique de la
repr\'esentation r\'eflexion (pour les groupes de type $A$, les deux
sont \'egaux).

Nous d\'efinissons maintenant deux entiers qui sont associ\'es \`a une
repr\'esentation de $W$ quelconque. La \emph{$a$-valeur} d'une
repr\'esentation est l'entier $i$ le plus petit tel que le
coefficient de $q^i$ soit non z\'ero dans le degr\'e fant\^ome.  La
\emph{$\ta$-valeur} est la m\^eme quantit\'e \`a l'\'egard du degr\'e
g\'en\'erique. (Nous suivons ici les notations de Carter \cite{CFG}~;
ailleurs dans la litt\'erature, ces entiers sont appel\'es la
$b$-valeur et la $a$-valeur, respectivement). D'apr\`es Lusztig, une
repr\'esentation est dite \emph{sp\'eciale} si sa $a$-valeur et sa
$\ta$-valeur sont \'egales. Les classes unipotentes sp\'eciales sont
les classes $\orb$ telles que $(\orb,1)$ correspond \`a une
repr\'esentation sp\'eciale via la correspondance de Springer.  Pour
toute classe unipotente $\orb$, nous dirons que $\nu^{-1}(\orb,1)$ est
\emph{la repr\'esentation de Springer de $\orb$}.

Le quotient de Lusztig de $A(\orb)$, not\'e $\Ab(\orb)$, est d\'efini
en termes des $a$-valeurs et des $\ta$-valeurs. Pour une classe $\orb$
donn\'ee, nous consid\'erons toutes les repr\'esentations de $W$ qui
correspondent \`a des couples $(\orb,\pi)$. En particulier, soit $E_0$
la repr\'esentation de $W$ associ\'ee \`a $(\orb,1)$. Soit $K \subset
A(\orb)$ l'intersection des noyaux des repr\'esentations $\pi$ telles
que $(\orb,\pi) = \nu(E)$ o\`u $E$ est une repr\'esentation de $W$
dont la $\ta$-valeur est \'egale \`a celle de $E_0$.  Le quotient de
Lusztig est le quotient de $A(\orb)$ par $K$ (\`a l'origine, Lusztig
fit cette d\'efinition seulement pour les classes sp\'eciales, mais
Sommers remarqua qu'elle est \'egalement valable pour les classes non
sp\'eciales).

\subsection{Cellules bilat\`eres}

Rappelons maintenant quelques faits sur les cellules bilat\`eres. Il
s'agit de certains sous-ensembles de l'ensemble des \'el\'ements de
$W$, d\'efinis via l'action de $W$ sur son alg\`ebre de Hecke. Les
\'el\'ements d'une cellule bilat\`ere correspondent \`a une base pour
un module (non irr\'eductible, en g\'en\'eral) de $W$, et on peut donc
parler d'une repr\'esentation qui ``appara\^\i t'' dans une cellule
bilat\`ere. En fait toute repr\'esentation de $W$ appara\^\i t dans une
unique cellule bilat\`ere et chaque cellule bilat\`ere contient une
unique repr\'esentation sp\'eciale.

Lusztig a trouv\'e une relation int\'eressante entre les cellules
bilat\`eres et les repr\'esentations de Springer, fond\'ee sur la
notion de \emph{pi\`ece sp\'eciale}. Une pi\`ece sp\'eciale est la
r\'eunion d'une classe unipotente sp\'eciale et toutes les classes
unipotentes dans son adh\'erence de Zariski qui n'appartiennent \`a
l'adh\'erence de Zariski d'aucune classe sp\'eciale plus petite. Donc
une pi\`ece sp\'eciale contient une classe sp\'eciale et un certain
nombre (\'eventuellement nul) de classes non sp\'eciales. Spaltenstein
a montr\'e que toute classe unipotente appartient \`a une unique
pi\`ece sp\'eciale, et donc que les pi\`eces sp\'eciales constituent
une partition de la vari\'et\'e unipotente. Le th\'eor\`eme suivant,
qui relie les pi\`eces sp\'eciales aux cellules bilat\`eres, est un
corollaire imm\'ediat d'un r\'esultat de Lusztig
{\cite[Theorem~0.2]{LNU}}.

\begin{thm}\label{thm:lusztig}
Soient $\orb_1$ et $\orb_2$ deux classes unipotentes, et soient $E_1$
et $E_2$ leurs repr\'esentations de Springer respectives. $\orb_1$ et
$\orb_2$ appartiennent \`a la m\^eme pi\`ece sp\'eciale si et
seulement si $E_1$ et $E_2$ appartiennent \`a la m\^eme cellule
bilat\`ere.
\end{thm}

Si l'une des repr\'esentations $E_i$ correspond \`a un couple
$(\orb_i,\pi_i) \in \N_G$ avec $\pi_i$ non triviale, l'analogue de
l'\'enonc\'e ci-dessus n'est pas vrai: il est possible, dans ce cas,
que $\orb_1$ et $\orb_2$ n'appartiennent pas \`a la m\^eme pi\`ece
sp\'eciale. N\'eanmoins, Geck et Malle ont obtenu le r\'esultat
suivant dans ce contexte.

\begin{prop}[{\cite[Proposition~2.2]{GMU}}]\label{prop:cell}
Soit $E$ une repr\'esentation sp\'eciale de $W$ qui correspond \`a la
classe $\orb$. Si $E_1$ est une repr\'esentation irr\'eductible dans
la cellule bilat\`ere de $E$, alors $E_1$ est associ\'ee par la
correspondance de Springer \`a un couple $(\orb_1,\pi_1)$ tel que
$\orb_1 \leq \orb$.
\end{prop}

\subsection{Familles de Lusztig}
\label{sect:familles}

Lusztig a d\'efini une partition de $\Irr(W)$ en \emph{familles} de la
mani\`ere suivante. Si $W=\{1\}$, il y a une famille unique, constitu\'ee
par la repr\'esentation triviale de $W$. Supposons maintenant que $W$
poss\`ede au moins deux \'el\'ements et que les familles ont d\'ej\`a
\'et\'e d\'efinies pour tous les sous-groupes paraboliques standard
propres de $W$. On dit alors deux
repr\'esentations irr\'eductibles $E$ et $E'$ de $W$ appartiennent \`a la
m\^eme famille s'il existe une suite
$E=E_0$, $E_1$, $\ldots$, $E_r=E'$ de repr\'esentations irr\'eductibles de
$W$, telles que, pour chaque $i$ ($0\le i\le r-1$), 
il existe un sous-groupe parabolique propre standard
$W_i$ de $W$ et des repr\'esentations irr\'eductibles
$M_i'$, $M_i''$ de $W_i$, appartenant \`a une m\^eme famille de $W_i$,
telles que
$$
\begin{cases}
\langle M_i',E_{i-1}\rangle_{W_i}\ne 0,\ta_{M_i'}=\ta_{E_{i-1}},\cr
\langle M_i'',E_i\rangle_{W_i}\ne 0, \ta_{M_i''}=\ta_{E_i}
\end{cases}\;
\text{ou}\;
\begin{cases}
\langle M_i',E_{i-1}\otimes\sgn\rangle_{W_i}\ne 0,\ta_{M_i'}=
\ta_{E_{i-1}\otimes\sgn},\cr
\langle M_i'',E_i\otimes\sgn\rangle_{W_i}\ne 0,
\ta_{M_i''}=\ta_{E_i\otimes\sgn},
\end{cases}
$$
o\`u $\sgn$ est la repr\'esentations signe du groupe $W$.  La
fonction $\ta$ est constante sur chaque famille. Barbasch et Vogan ont
montr\'e que deux repr\'esentations irr\'eductibles de $W$
appartiennent \`a une m\^eme famille si et seulement si elles
apparaissent dans une m\^eme cellule bilat\`ere,
\cite[Theorem~2.29]{BV2} (voir aussi \cite{BV1}). Nous dirons que cette 
cellule et cette famille \emph{se correspondent}.
Une cellule bilat\`ere sera dite \emph{cuspidale} si la famille qui lui 
correspond est cuspidale au sens de \cite[(8.1)]{Lbook}.

\subsection{Induction}

Pour l'\'etude des supports des faisceaux caract\`eres, nous voudrions
conna\^\i tre des propri\'et\'es des repr\'esentations de $W$ induites
d'une certaine classe de sous-groupes de $W$. Il est commode
d'introduire en m\^eme temps une certaine classe analogue de
sous-groupes de $G$.  Dans la terminologie introduite par
Sommers, un sous-groupe de $G$ est appel\'e un \emph{pseudo
sous-groupe de Levi} s'il est le centralisateur d'un \'el\'ement
semi-simple. Un sous-groupe de $W$ qui est le groupe de Weyl d'un
pseudo sous-groupe de Levi sera appel\'e un \emph{pseudo sous-groupe
parabolique}. Une autre description est obtenue comme suit. Supposons
choisi un ensemble $\Pi$ de racines positives simples. Soit $\Pi_0 =
\Pi \cup \{\alpha_0\}$, o\`u $\alpha_0$ est le n\'egatif de la racine
la plus haute. Rappelons qu'un sous-groupe de $G$ (resp.~de $W$) est
un sous-groupe de Levi (resp.~parabolique) s'il correspond \`a un
sous-syst\`eme de racines engendr\'e par un sous-ensemble de $\Pi$.
Un sous-groupe de $G$ (resp.~de $W$) est un pseudo sous-groupe de Levi
(resp.~pseudo sous-groupe parabolique) s'il correspond \`a un
sous-syst\`eme de racines engendr\'e par un sous-ensemble propre de
$\Pi_0$.

Une classe de sous-groupes de $W$ qui nous sera tr\`es importante est
celle fournie par l'ensemble des pseudo sous-groupes de Levi du groupe
dual $G^*$ de $G$. Puisque le groupe de Weyl de $G$ et celui de $G^*$
sont canoniquement isomorphes, nous pouvons consid\'erer les pseudo
sous-groupes paraboliques du groupe de Weyl de $G^*$ comme
sous-groupes de $W$ lui-m\^eme.  Par abus de langage, nous appelerons
de tels sous-groupes de $W$ \emph{pseudo sous-groupes paraboliques du
dual de $W$}. De plus, la correspondance entre les sous-groupes de
Levi de $G$ et ceux de $G^*$ montre que tout sous-groupe parabolique
de $W$ est aussi un pseudo sous-groupe parabolique de son dual.

Si $W'$ est un pseudo sous-groupe parabolique du dual de $W$ et $E$
est une repr\'esentation irr\'eductible de $W'$, il y a deux
sous-modules de $\Ind_{W'}^W E$ auxquels il faut pr\^eter
attention. Notons $j_{W'}^W E$ la somme des composantes
irr\'eductibles de $\Ind_{W'}^W E$ qui ont la m\^eme $a$-valeur que
$E$. La somme de ceux qui ont la m\^eme $\ta$-valeur que $E$ sera
not\'ee $J_{W'}^W E$. Sous certaines hypoth\`eses, qui sont toujours
satisfaites si, par exemple, $E$ est sp\'eciale, le module $j_{W'}^W
E$ est lui aussi irr\'eductible. Plus pr\'ecis\'ement, si $E$ est
sp\'eciale, ce module, appel\'e \emph{l'induction tronqu\'ee} de $E$,
est toujours une repr\'esentation de Springer d'une certaine classe
unipotente. De plus, si $W'$ est un sous-groupe parabolique, alors
l'induction tronqu\'ee de $E$ est aussi sp\'eciale.

L'induction pour les repr\'esentations des groupes de Weyl est
intimement li\'ee \`a certains autres types d'induction. L'induction
pour les classes unipotentes est pr\'ecis\'ement l'induction
tronqu\'ee des repr\'esentations de Springer correspondantes. De
surcro\^it, pour une cellule bilat\`ere d'un pseudo sous-groupe
parabolique du dual donn\'ee, on peut consid\'erer la cellule qui
contient la $j$-induction de son unique repr\'esentation
sp\'eciale. Appelons cette cellule la cellule \emph{induite} de la
premi\`ere. (Cette terminologie est en accord avec la th\'eorie
d'induction plus sophistiqu\'ee introduite par Xi dans \cite{XIC}).
La repr\'esentation $J_{W'}^W E$ est telle que toutes ses composantes
irr\'eductibles appartiennent \`a une unique cellule bilat\`ere, qui
est la cellule induite de celle \`a laquelle appartient $E$.

L'ensemble des cellules bilat\`eres est muni d'un ordre partiel
naturel provenant de leur d\'efinition en termes de l'action de $W$
sur son alg\`ebre de Hecke. D'autre part, cet ensemble est en
bijection avec l'ensemble des classes unipotentes sp\'eciales, qui
h\'erite un ordre partiel de l'ensemble de toutes les classes
unipotentes. Barbasch et Vogan ont d\'emontr\'e que ces deux ordres
co\"\i ncident \cite[Proposition~3.23]{BV3}. Cette \'equivalence
permet de voir que l'induction des cellules bilat\`eres est une
application croissante, en la comparant avec l'induction des classes
unipotentes. D'apr\`es Spaltenstein, cette derni\`ere op\'eration est
croissante. Soient $\orb_1$ et $\orb_2$ deux classes induites telles
que $\orb_1 \leq \orb_2$, et soient $\orb'_1$ et $\orb'_2$ les classes
sp\'eciales dans leurs pi\`eces sp\'eciales respectives. Pour voir que
l'induction des cellules est croissante, il faut d\'emontrer que
$\orb'_1 \leq \orb'_2$, laquelle in\'egalit\'e est impliqu\'ee par
certaines propri\'et\'es de l'application de dualit\'e $d$ de
Spaltenstein. Cette application est d\'ecroissante, donc $d^2(\orb_1)
\leq d^2(\orb_2)$. De plus, $d^2(\orb)$ est toujours la plus petite classe
unipotente sp\'eciale dont l'adh\'erence de Zariski contient $\orb$.
Autrement dit, $d^2(\orb)$ est l'unique classe sp\'eciale dans la pi\`ece 
sp\'eciale de $\orb$. On d\'eduit que $\orb'_1 \leq \orb'_2$.

\section{Repr\'esentations bien support\'ees}
\label{sect:biensupp}

Via la correspondance de Springer, nous pouvons introduire la notion de 
\emph{support} d'une repr\'esentation de $W$. Soit $E$ une
repr\'esentation de $W$. Si $E$ irr\'eductible et $\nu(E) =
(\orb,\pi)$, le support de $E$ est la classe unipotente $\orb$. Si $E$
n'est pas irr\'eductible, son support est la r\'eunion des supports de
ses composantes irr\'eductibles. Les r\'esultats de cette section sont
fond\'es sur la notion suivante.

\begin{defn}\label{defn:b-s}
Une repr\'esentation $E$ de $W$ est dite \emph{bien support\'ee} si
son support $\supp E$ v\'erifie les deux conditions suivantes:
\begin{enumerate}
\item Il existe une unique classe unipotente $\orb_0$ telle que la
repr\'esentation de Springer $E_0$ de $\orb_0$ intervienne dans $E$, et
$\orb_0 \subset \supp E \subset
\overline{\orb_0}$. \label{enum:bien-supp1}
\item Tous les termes de $E$ appartiennent \`a des cellules bilat\`eres
$\bc'$ qui satisfont $\bc' \leq \bc_0$, o\`u $\bc_0$ d\'esigne la
cellule qui contient $E_0$. \label{enum:bien-supp2}
\end{enumerate}
De plus, $E$ est dite \emph{sp\'ecialement bien support\'ee} si la
classe $\orb_0$ est sp\'eciale.
\end{defn}

Une cons\'equence particuli\`ere de cette d\'efinition sera importante
dans la suite. Pour le lemme suivant, nous conservons les notations de
la d\'efinition pr\'ec\'edente.

\begin{lem}\label{lem:ab-defini}
Soit $E_1$ une composante irr\'eductible de $E$, et supposons que
$\nu(E) = (\orb,\pi)$.  Si $\orb$ appartient \`a la pi\`ece sp\'eciale
de $\orb_0$, alors $\pi$ est d\'efinie sur $\Ab(\orb)$.
\end{lem}
\begin{proof}
Gr\^ace \`a la condition~(\ref{enum:bien-supp2}), nous savons que
$\nu^{-1}(\orb,\pi)$ doit appartenir \`a $\bc_0$, car selon la
Proposition~\ref{prop:cell}, les repr\'esentations dans une cellule
plus petite ne peuvent pas \^etre associ\'ees \`a une classe dans
la pi\`ece sp\'eciale de $\orb_0$. Par cons\'equent, la $\ta$-valeur
de $\nu^{-1}(\orb,\pi)$ est \'egale \`a celle de $E_0$. Le
Th\'eor\`eme~\ref{thm:lusztig} nous dit que cette derni\`ere est \'egale
\`a la $\ta$-valeur de la repr\'esentation de Springer de
$\orb$. On d\'eduit que le noyau $K$ de l'application $A(\orb) \to
\Ab(\orb)$ est contenu dans le noyau de $\pi$, et donc que $\pi$ est
d\'efinie sur $\Ab(\orb)$.
\end{proof}

La proposition suivante relie le support d'une repr\'esentation
induite et les ordres partiels sur les classes et sur les
cellules. 

\begin{prop}\label{prop:ind-supp}
Soit $L^*$ un pseudo sous-groupe de Levi de $G^*$, et soit $W'$ son
groupe de Weyl. Soit $E_1$ une repr\'esentation irr\'eductible de $W'$
qui appartient \`a une cellule bilat\`ere $\bc_1$ correspondant \`a la
classe unipotente sp\'eciale $\orb_1$. Soit $\bc$ et $\orb$ les
induites de $\bc_1$ et $\orb_1$, respectivement, et soit $E =
\Ind_{W'}^W E_1$.  Alors:
\begin{enumerate}
\item $\supp E \subset \overline{\orb}$. \label{enum:ind-supp1}
\item Toute composante irr\'eductible de $E$ appartient \`a une cellule
bilat\`ere $\bc'$ telle que $\bc' \leq \bc$. \label{enum:ind-supp2}
\end{enumerate}
\end{prop}

La plupart de cette section est consacr\'ee \`a la preuve de cette
proposition. Nous l'\'etablirons pour les groupes classiques par des
m\'ethodes combinatoires, et pour les groupes exceptionnels par des
calculs explicites. Avant de commencer ce projet, cependant, nous
verrons la cons\'equence la plus importante de cette proposition.

\begin{thm}\label{thm:b-s}
Soit $L^*$ un pseudo sous-groupe de Levi de $G^*$, et soit $W'$ le
groupe de Weyl de $L^*$. Soit $E_1$ une repr\'esentation de $W'$, et
soit $E = \Ind_{W'}^W E_1$. Si $E_1$ est sp\'ecialement bien
support\'ee, alors $E$ est bien support\'ee. De plus, lorsque $L^*$ est
un sous-groupe de Levi, la repr\'esentation $E$ est sp\'ecialement
bien support\'ee.
\end{thm}
Il est \`a noter que $E_1$ n'est pas suppos\'ee irr\'eductible.
\begin{proof}
Soit $\orb_1$ une classe maximale dans le support de $E_1$, et soit
$\orb_0$ la classe induite de $\orb_1$. La
Proposition~\ref{prop:ind-supp} nous dit que $\supp E \subset
\overline{\orb_0}$. De plus, puisque la repr\'esentation de Springer
de $\orb_1$ intervient dans $E_1$, celle de $\orb_0$ intervient dans
$E$, et par cons\'equent $\orb_0 \subset \supp E$. La
condition~(\ref{enum:bien-supp1}) d'\^etre bien support\'ee est donc
satisfaite.

Soit $\bc_1$ la cellule bilat\`ere qui correspond \`a $\orb_1$, et
soit $\bc_0$ sa cellule induite dans
$W$. L'\'enonc\'e~(\ref{enum:ind-supp2}) de la
Proposition~\ref{prop:ind-supp}, combin\'e avec le fait que
l'induction des cellules bilat\`eres respecte l'ordre partiel,
implique la condition~\ref{enum:bien-supp2}.

La repr\'esentation $E$ est donc bien support\'ee. Si $L$ est un
sous-groupe de Levi, nous savons que la classe induite $\orb_0$ est
sp\'eciale, et donc sous cette hypoth\`ese $E$ est sp\'ecialement bien
support\'ee.
\end{proof}

Consid\'erons maintenant les aspects d'un groupe classique qu'il faut
comprendre pour prouver la Proposition~\ref{prop:ind-supp}.  Nous
employerons certains objets combinatoires introduits par Lusztig pour
param\'etrer les repr\'esentations de $W$ et les \'el\'ements de
$\N_G$. Il s'agit des symboles \cite{LIR} et des $u$-symboles
\cite{LIC}, respectivement.  Pour l'instant, nous restreignons notre
attention aux $u$-symboles qui correspondent aux \'el\'ements de
$\N_G$ apparaissant dans la correspondance de Springer (non
g\'en\'eralis\'ee): ce sont les $u$-symboles ``de d\'efaut ${}\leq
1$''. Les symboles et les $u$-symboles de d\'efaut $\leq 1$ sont
certaines paires de listes d'entiers, de la forme
\[
\underset{\text{types $B$ et $C$}}{\begin{symbole}
a_1&  &a_2&  &\cdots&  &a_{m+1}\\
&b_1&  &\cdots&  &b_m
\end{symbole}}
\qquad\text{ou}\qquad
\underset{\text{type $D$}}{\begin{symbole}
a_1&\quad &a_2&\quad &\cdots&\quad &a_m\\
b_1&      &b_2&      &\cdots&      &b_m
\end{symbole}}
\]
o\`u $a_1 < a_2 < \cdots < a_{m+1}$, et $b_1 < \cdots < b_m$.  Deux
symboles ou $u$-symboles dans lequels les m\^emes entiers ont lieu
avec les m\^emes multiplicit\'es sont dits \emph{similaires}.  Un
symbole ou un $u$-symbole est dit \emph{distingu\'e} si
\[
a_1 \leq b_1 \leq a_2 \leq b_2 \leq \cdots \leq a_m \leq b_m \leq
a_{m+1}.
\]
Il est \'evident que toute classe de similitude contient un unique
\'el\'ement distingu\'e. Un symbole correspond \`a une
repr\'esentation sp\'eciale si et seulement s'il est distingu\'e.  Un
$u$-symbole correspond \`a un couple $(\orb,1)$ si et seulement s'il
est distingu\'e.

Appelons le nombre $m$ la \emph{longueur} du symbole ou du
$u$-symbole. Il y a une relation d'\'equivalence qui permet, pour tout
symbole ou $u$-symbole de longueur $m$, de trouver un symbole ou
$u$-symbole \'equivalent de longueur $m+1$.

Dans le type $B_n$, soit $N = 2n+1$; dans $C_n$ ou $D_n$, soit $N =
2n$. Les classes unipotentes sont ainsi en correspondance avec un
certain sous-ensemble de l'ensemble des partitions de $N$.  Si
$\lambda = (\lambda_1 \leq \lambda_2 \leq \cdots \leq \lambda_k)$ est
une partition de $N$, soit $\sigma_i(\lambda)$ la somme des $i$ plus
grandes parties de $\lambda$: c'est-\`a-dire,
\[
\sigma_i(\lambda) = \sum_{j=k-i+1}^k \lambda_j.
\]
Dans l'ordre partiel usuel sur les partitions, on dit que $\lambda
\leq \lambda'$ si $\sigma_i(\lambda) \leq \sigma_i(\lambda')$ pour tout $i$.
Si $\lambda=(\lambda_1 \leq \lambda_2 \leq \cdots \leq \lambda_k)$ et
$\lambda'=(\lambda_1 \leq \lambda_2 \leq \cdots \leq \lambda_k)$ sont
des partitions de $n$ et $n'$ respectivement, nous noterons $\lambda
\vee \lambda'$ la partition de $n+n'$ d\'efinie par
$(\lambda\vee\lambda')_i = \lambda_i+\lambda'_i$.
Nous dirons que $\lambda\vee\lambda'$ est la partion \emph{jointe} de
$\lambda$ et $\lambda'$.

Nous d\'ecrivons maintenant le processus pour obtenir un symbole ou un
$u$-symbole \`a partir d'une partition de $N$. Le nombre de parties de
$\lambda$ est pair dans le type $D$ et impair dans le type
$B$. Supposons ce nombre pair dans le type $C$, quitte \`a ajouter une
partie \'egale \`a $0$ si n\'ecessaire. Soit $m$ tel que le nombre de
parties de $\lambda$ est $2m$ (types $C$ et $D$) ou $2m+1$ (type $B$),
et soit $\bar\lambda$ la partition d\'efinie par la formule
$\bar\lambda_i = \lambda_i + (i-1)$: c'est une partition de $N +
m(2m-1)$ (types $C$ et $D$) ou $N + m(2m+1)$ (type $B$).

Soit $\eta^* = (\eta^*_1 < \cdots < \eta^*_m)$ la liste des parties
paires de $\bar\lambda$ (on peut montrer par r\'ecurrence que $\eta^*$
doit avoir $m$ parties). Soit $\xi^*$ la liste des parties impaires de
$\bar\lambda$, avec un ``$1$'' suppl\'ementaire ajout\'e au d\'ebut de
la liste dans le type $C$. D\'efinissions deux partitions $\eta$ et
$\xi$ telles que $2\eta_i = \eta^*_i$ et $2\xi_i+1 = \xi^*_i$. Alors
$\xi$ a $m+1$ parties dans les types $B$ et $C$, et $m$ parties dans
le type $D$. Soient $\bar\xi$ et $\bar\eta$ les listes d'entiers
d\'efinies par les formules suivants.
\begin{align*}
\text{Type $B$~:}&& 
\binom{\bar\xi}{\bar\eta} &=
\begin{symbole}
\xi_1& &\xi_2+1& &\cdots& &\xi_m+(m-1)& &\xi_{m+1}+m \\
&\eta_1& &\eta_2+1& &\cdots& &\eta_m+(m-1)
\end{symbole}
\\
\text{Type $C$~:}&& 
\binom{\bar\xi}{\bar\eta} &=
\begin{symbole}
\xi_1& &\xi_2+2& &\cdots& &\xi_m+m& &\xi_{m+1}+(m+1) \\
&\eta_1+1& &\eta_2+2& &\cdots& &\eta_m+m
\end{symbole}
\\
\text{Type $D$~:}&& 
\binom{\bar\xi}{\bar\eta} &=
\begin{symbole}
\xi_1&\quad &\xi_2+1&\quad &\cdots&\quad &\xi_m+(m-1) \\
\eta_1&\quad &\eta_2+1&\quad &\cdots&\quad &\eta_m+(m-1) \\
\end{symbole}
\end{align*}
Si $\lambda$ est la partition qui correspond \`a la classe unipotente
$\orb$, alors $\binom{\bar\xi}{\bar\eta}$ est le $u$-symbole qui
correspond \`a $(\orb,1)$, et $\binom{\xi}{\eta}$ est le symbole qui
correspond \`a la repr\'esentation du groupe de Weyl associ\'ee \`a
$(\orb,1)$.

Il peut arriver qu'une comparaison de deux symboles ou $u$-symboles
implique une relation d'ordre entre les classes unipotentes
associ\'ees. En particulier, remarquons que si $\lambda$ et $\lambda'$
sont deux partitions de $N$ correspondant \`a deux classes
unipotentes, alors $\lambda \leq \lambda'$ si et seulement si
$\bar\lambda \leq \bar\lambda'$. De plus, ces derni\`eres partitions
sont faciles \`a calculer \`a partir du symbole correspondant.  Cette
observation est utilis\'ee au cours de la preuve du lemme suivant.

\begin{lem}\label{lem:ordre}
Soient $\binom{\bar\xi}{\bar\eta}$ et $\binom{\bar\xi'}{\bar\eta'}$
les $u$-symboles de la m\^eme longueur, et soient $\lambda$ et
$\lambda'$ les partitions param\'etrant les classes unipotentes
auxquelles ils sont associ\'es. Si $\bar\xi \cup \bar\eta \leq
\bar\xi' \cup \bar\eta'$, alors $\lambda \leq \lambda'$, avec
\'egalit\'e si et seulement si $\binom{\bar\xi}{\bar\eta} =
\binom{\bar\xi'}{\bar\eta'}$.
\end{lem}
\begin{proof}
Nous pouvons supposer les $u$-symboles distingu\'es, car les
partitions $\bar\xi \cup \bar\eta$ et $\bar\xi' \cup \bar\eta'$ ne
d\'ependent que de la classe de similitude des $u$-symboles. Ces
$u$-symboles proviennent alors de deux symboles $\binom{\xi}{\eta}$ et
$\binom{\xi'}{\eta'}$.
D'apr\`es les commentaires ci-dessus, il suffit d'\'etablir que
$\bar\lambda \leq \bar\lambda'$. Supposons que les partitions $\bar\xi
\cup \bar\eta$ et $\bar\xi' \cup \bar\eta'$ ne diff\`erent qu'en deux
parties (pour tout couple de $u$-symboles, il est possible de trouver
une suite de $u$-symboles interm\'ediaires tels que deux $u$-symboles
cons\'ecutifs ne diff\`erent jamais en plus de deux parties). On en
d\'eduit que les symboles $\binom{\xi}{\eta}$ et $\binom{\xi'}{\eta'}$
diff\`erent aussi en deux parties au plus.

Notons les partitions $\xi \cup \eta$ et $\xi' \cup \eta'$
respectivement $\mu$ et $\mu'$. Il existe donc deux indices $i_0 >
i_1$ telles que $\mu_i = \mu'_i$ si $i \neq i_0,i_1$, mais $\mu_{i_0}
= \mu'_{i_0} - C$ et $\mu_{i_1} = \mu'_{i_1} + C$, o\`u $C$ est un
certain entier strictement positif. Cette description nous fournit une
description analogue de la relation entre $\bar\lambda$ et
$\bar\lambda'$. Soit $\bar\lambda'_{j_0}$ la partie de $\bar\lambda'$
qui correspond \`a $\mu'_{i_0}$: c'est-\`a-dire, ou $\bar\lambda'_{j_0} =
2\mu'_{i_0}+1$ ou $\bar\lambda'_{j_0} = 2\mu'_{i_0}$, suivant que
$\mu'_{i_0}$ est membre de $\xi'$ ou de $\eta'$. Soit
$\bar\lambda'_{j_1}$ la partie analogue correspondant \`a $\mu'_{i_1}$
(il n'est pas forc\'ement vrai que $i_0 = j_0$ ou $i_1 = j_1$). Alors
$\bar\lambda$ est la partition obtenue \`a partir de $\bar\lambda'$ en
rempla\c cant $\bar\lambda'_{j_0}$ et $\bar\lambda'_{j_1}$ par
$\bar\lambda'_{j_0} - 2C$ et $\bar\lambda'_{j_1} + 2C$,
respectivement. Remarquons, cependant, que ce n'est pas dire que
$\bar\lambda_{j_0} = \bar\lambda'_{j_0} - 2C$. Apr\`es d'avoir
remplac\'e ces deux parties de $\bar\lambda'$, il est possible que les
parties de celle-ci ne soient plus en l'ordre d\'ecroissant. Donc il
n'est pas tout de suite \'evident que $\bar\lambda \leq \bar\lambda'$.

Nous savons que $\mu_{i_0}-C \geq \mu_{i_1} +C$. Si cette
in\'egalit\'e est stricte, alors on peut conclure que
$\bar\lambda'_{j_0}-2C \geq \bar\lambda'_{j_1}-2C$. D'autre part, si
la premi\`ere in\'egalit\'e est en r\'ealit\'e une \'egalit\'e, alors
il est possible que $\bar\lambda'_{j_0}-2C = \bar\lambda_{j-1}-2C-1$,
dans le cas o\`u $\mu'_{i_0}$ fait partie de $\bar\eta'$ et
$\mu'_{i_1}$ fait partie de $\bar\xi'$. N\'eanmoins, ces deux parties,
consid\'er\'ees toutes seules comme partitions de $\bar\lambda'_{j_0}
+ \bar\lambda'_{j_1}$, satisfont toujours
\[
[\bar\lambda'_{j_1} < \bar\lambda'_{j_0}] <
[\bar\lambda'_{j_1}+2C,\bar\lambda'_{j_0}-2C]
\]
(l'in\'egalit\'e est stricte parce que $C$ est suppos\'e
strictement positif). Toutes les autres parties de $\bar\lambda$ et de
$\bar\lambda'$ sont \'egales, donc nous pouvons conclure que
$\bar\lambda < \bar\lambda'$.
\end{proof}

\begin{lem}\label{lem:comp}
Soient $(\alpha,\beta)$ et $(\alpha',\beta')$ deux paires de
partitions, vues comme param\'etrant des repr\'esentations de $W$,
telles que $\alpha \leq \alpha'$ et $\beta \leq \beta'$.
\begin{enumerate}
\item Soient $(\orb,\pi)$ et $(\orb',\pi')$ les \'el\'ements de $\N_G$
auxquels elles sont associ\'ees par la correspondance de
Springer. Alors $\orb \leq \orb'$, avec \'egalit\'e seulement si
$\alpha = \alpha'$ et $\beta = \beta'$. \label{lem:comp1}
\item Soient $\orb_1$ et $\orb'_1$ les classes sp\'eciales
auxquelles correspondent les cellules bilat\`eres dont
$(\alpha,\beta)$ et $(\alpha',\beta')$ font partie, respectivement.
Alors $\orb_1 \leq \orb'_1$, avec \'egalit\'e seulement si $\alpha =
\alpha'$ et $\beta = \beta'$. \label{lem:comp2}
\end{enumerate}
\end{lem}
\begin{proof}
Soient $\binom{\xi}{\eta}$, $\binom{\xi'}{\eta'}$ et
$\binom{\bar\xi}{\bar\eta}$, $\binom{\bar\xi'}{\bar\eta'}$ les
symboles et les $u$-symboles qui correspondent \`a $(\alpha,\beta)$ et
$(\alpha',\beta')$, respectivement. Supposons que tous ces symboles et
$u$-symboles ont m\^eme longueur.

Associons \`a $(\alpha,\beta)$ et $(\alpha',\beta')$ deux $u$-symboles
$\binom{\bar\xi}{\bar\eta}$ et $\binom{\bar\xi'}{\bar\eta'}$ de la
m\^eme longueur. Gr\^ace aux hypoth\`eses sur $(\alpha,\beta)$ et
$(\alpha',\beta')$, il est clair que $\bar\xi \cup \bar\eta \leq
\bar\xi' \cup \bar\eta'$. L'\'enonc\'e (\ref{lem:comp1}) est alors
cons\'equence imm\'ediate du Lemme~\ref{lem:ordre}.

Quant \`a l'\'enonc\'e (\ref{lem:comp2}), remarquons d'abord qu'il
n'\'equivaut pas \`a l'\'enonc\'e que l'unique classe sp\'eciale dans
la pi\`ece sp\'eciale de $\orb$ est inf\'erieure \`a celle dans la
pi\`ece sp\'eciale de $\orb'$. La r\'eciproque de la
Proposition~\ref{prop:cell} est fausse: il est possible que
$\orb_1$ soit strictement plus grande que la classe sp\'eciale dans la
pi\`ece sp\'eciale de $\orb$.

La partie~(\ref{lem:comp2}) d\'ecoule d'un bref calcul en termes des
symboles. Soient $\lambda$ et $\lambda'$ les partitions des classes
$\orb_1$ et $\orb'_1$, et soient $\bar\lambda$ et $\bar\lambda'$
d\'efinies d'apr\`es la discussion qui pr\'ec\`ede le
Lemme~\ref{lem:ordre}. On a $\bar\lambda = \xi \cup \eta$ et $\xi'
\cup \eta'$. Il est \'evident, sous les hypoth\`eses du lemme, que
$\xi \leq \xi'$ et $\eta \leq \eta'$, et donc $\bar\lambda \leq
\bar\lambda'$. Il est \'egalement clair que $\bar\lambda =
\bar\lambda'$ seulement si $\xi = \xi'$ et $\eta = \eta'$. 
\end{proof}

\begin{lem}\label{lem:comp-A}
Soit $n_1 = n/2$ si $n$ est pair, et $n_1 = (n+1)/2$ si $n$ est
impair. Soit $n_2 = n- n_1$. Soient $(\alpha,\beta) =
([1^{n_1}],[1^{n_2}])$ et $(\alpha',\beta') =
([1^{n_1+k}],[1^{n_2-k}])$, o\`u $-n_1 \leq k \leq n_2$. Alors les
conclusions (\ref{lem:comp1}) et (\ref{lem:comp2}) du
Lemme~\ref{lem:comp} sont vraies.
\end{lem}
\begin{proof}
Nous allons v\'erifier l'\'enonc\'e~(\ref{lem:comp1}) du lemme
pr\'ec\'edent par un calcul explicite des $u$-symboles
appropri\'es. Supposons pour l'instant que $k$ est positif, et que le
groupe $W$ est de type $B$. Alors
\[
\binom{\bar\xi}{\bar\eta'} =
\begin{symbole}
0&\quad &2&\quad &\cdots \quad l\quad l+3\quad \cdots\quad 2m-1& &2m+1\\
&0& &2& \cdots \quad p\quad p+3\quad \cdots& 2m-1\\
\end{symbole},
\]
o\`u $l = 2(m-n_1-k)$ et $p = 2(m-n_2+k)$. La partition $\bar\xi' \cup
\bar\eta'$ a donc la forme
\begin{multline*}
[0^2,2^2,\ldots,l^2,l+2,l+3,\ldots,p-1,p,\\
(p+2)^2,(p+4)^2,\ldots,2(m-2)+1^2,2(m-1)+1^2,2m+1].
\end{multline*}
Remarquons que $p = l + 2(2k+n_1-n_2)$. Il est \'evident, \`a partir
de cette description, que la partition ci-dessus atteint la valeur la
plus haute possible (\`a l'\'egard de l'ordre partiel sur les
partitions) quand $k = 0$. Nous appliquons ensuite le
Lemme~\ref{lem:ordre} pour obtenir le r\'esultat requis.

Pour \'etablir l'\'enonc\'e~(\ref{lem:comp2}), on r\'ep\`ete le calcul
pr\'ec\'edent en termes des symboles au lieu des $u$-symboles.  De
plus, il est \'egalement facile de traiter le cas o\`u $k < 0$, ainsi
que ceux de type $C$ et $D$.
\end{proof}

Nous retournons maintenant \`a la preuve de la
Proposition~\ref{prop:ind-supp}. Comme nous avons d\'ej\`a remarqu\'e,
la d\'emonstration dans le cas d'un groupe exceptionnel ne consiste
qu'en des calculs explicites des repr\'esentation induites. Cependant,
il n'est pas n\'ecessaire de calculer l'induite de toute
repr\'esentation de toute pseudo sous-groupe parabolique du groupe
dual. La premi\`ere moiti\'e de l'argument ci-dessous, qui est de
toute fa\c con n\'ecessaire pour les groupes classiques, permet aussi
de diminuer le nombre de calculs explicites qu'il faut faire pour les
groupes exceptionnels.

\begin{proof}[D\'emonstration de la Proposition~\ref{prop:ind-supp}]
Nous commen\c cons par d\'emontrer qu'il suffit de consid\'erer les
sous-groupes de Levi et pseudo sous-groupes de Levi maximaux de $G^*$
(nous verrons au cours de cette d\'emonstration la raison pour
laquelle il ne suffirait pas de consid\'erer les seuls pseudo
sous-groupes de Levi maximaux). Si $L^*$ n'est pas un tel sous-groupe
maximal, soit $M^*$ un pseudo sous-groupe de Levi de $G^*$ tel que
$L^*$ soit un sous-groupe de Levi de $M^*$, et soit $W''$ son groupe
de Weyl. Soit $M$ le groupe dual de $M^*$. (Il est probable que $M$ ne
soit isomorphe \`a aucun sous-groupe de $G$). Supposons la proposition
d\'ej\`a \'etablie pour $L^*$ en tant qu'un pseudo sous-groupe de Levi
du dual de $M$, ainsi que pour $M^*$ comme pseudo sous-groupe de Levi
du dual de $G$. \'Ecrivons la repr\'esentation induite de $E_1$ comme
somme d'un nombre fini de repr\'esentations irr\'eductibles de $W''$:
\[
\Ind_{W'}^{W''} E_1 = \sum m_i F_i,
\]
o\`u les $m_i$ sont des entiers positifs. Soit $\bd_i$ la cellule
bilat\`ere de $W''$ \`a laquelle appartient $F_i$, et soit $\orb'_i$
la classe unipotente sp\'eciale de $M^*$ \`a laquelle correspond
$\bd_i$. En particulier, supposons les \'etiquettes $i$ choisies de
telle sorte que $\bd_1$ soit la cellule induite de $\bc_1$. Il est
important de remarquer ici que $\orb'_1$ est la classe induite de
$\orb_1$, puisque $L^*$ est un sous-groupe de Levi de $M^*$. On sait
que $\Ind_{M^*}^G \orb'_1 = \orb$, ainsi que $\Ind_{W''}^W \bd_1 =
\bc$.

Pour tout $i$, nous savons que $\supp \Ind_{W''}^W F_i \subset
\overline{\Ind_{M^*}^G \orb'_i} \subset \overline{\orb}$, o\`u la
deuxi\`eme inclusion est cons\'equence du fait que l'induction est une
application croissante.  L'\'enonc\'e~(\ref{enum:ind-supp1}) est ainsi
\'etabli. Les cellules $\bc'$ auxquelles appartiennent les termes de
$\Ind_{W''}^W F_i$ satisfont $\bc' \leq \Ind_{W''}^W \bd_i$, mais
puisque l'induction des cellules est aussi croissante, le fait que
$\bd_i \leq \bd_1$ implique que $\Ind_{W''}^W \bd_i \leq
\bc$. L'\'enonc\'e~(\ref{enum:ind-supp2}) est donc \'egalement \'etabli.

Dor\'enavant, nous supposons que $L^*$ est soit un sous-groupe de Levi
maximal soit un pseudo sous-groupe de Levi maximal. Si $G$ est simple
et de type exceptionnel, il faut simplement calculer les induites de
toutes les repr\'esentations irr\'eductibles de tous les tels
sous-groupes maximaux. Les auteurs ont effectu\'e ces calculs \`a
l'aide du logiciel \textsf{CHEVIE} \cite{CHEVIE}. Nous ne donnons pas
de d\'etails de ces calculs.

Consid\'erons maintenant le cas o\`u $G$ est simple et de type
classique. Les sous-groupes de Levi et les pseudo sous-groupes de
Levi maximaux de $G^*$ ont les formes suivantes:
\begin{center}
\begin{tabular}{lcc}
           & sous-groupes de Levi     & pseudo sous-groupes de Levi \\
Type $B$: & $A_{k-1} \times C_{n-k}$ & $C_k \times C_{n-k}$ \\
Type $C$: & $A_{k-1} \times B_{n-k}$ & $D_k \times B_{n-k}$ \\
Type $D$: & $A_{k-1} \times D_{n-k}$ & $D_k \times D_{n-k}$
\end{tabular}
\end{center}
Nous pouvons en fait faire une r\'eduction suppl\'ementaire: dans le
type $B$, par exemple, on peut faire l'induction d'un sous-groupe de
Levi maximal en deux \'etapes, d'abord de $A_{k-1} \times C_{n-k}$ \`a
$C_k \times C_{n-k}$, et puis de ce dernier pseudo sous-groupe de Levi
\`a $B_n$. Il suffit donc de traiter les sous-groupes de Levi maximaux
de la forme $A_{n-1}$. De plus, il est seulement n\'ecessaire de
consid\'erer les repr\'esentations de $W'$ qui ne r\'esultent pas de
l'induction tronqu\'ee d'une repr\'esentation d'un sous-groupe
parabolique propre.

Pour le groupe de Weyl $W'$ de type $A_{n-1}$, la seule
repr\'esentation qui n'est pas induite de cette fa\c con est la
repr\'esentation signe, laquelle correspond \`a la partition
$[1^n]$. Nous savons que la multiplicit\'e de la repr\'esentation
correspondant \`a $(\alpha,\beta)$ dans $\Ind_{W'}^W [1^n]$ est
donn\'ee par le coefficient de Littlewood-Richardson
$c^{[1^n]}_{\alpha\beta}$. Ce dernier s'annule si $(\alpha,\beta)$ n'est
pas de la forme $([1^p],[1^{n-p}])$. Les deux conclusions du
Lemme~\ref{lem:comp-A}, traduites dans le langage des classes et des
cellules, deviennent les deux parties de la pr\'esente proposition.

Pour les pseudo sous-groupes de Levi maximaux, il n'est pas avantageux
d'exclure les repr\'esentations induites de notre discussion.  Nous
employons le Lemme~\ref{lem:ind-part} \'enonc\'e ci-dessous. Ce
lemme dit pr\'ecis\'ement ce qui est n\'ecessaire pour appliquer le
Lemme~\ref{lem:comp}, lequel implique ensuite la pr\'esente
proposition.
\end{proof}

Le m\^eme argument que celui utilis\'e dans la preuve de
\cite[Proposition 4.1]{Au} permet de d\'eduire le lemme suivant des
assertions VIII.3.(2), 4.(3), 5.(1) et 5.(2) de \cite{W}.

\begin{lem}\label{lem:ind-part}
Soit $L^*$ un pseudo sous-groupe de Levi de $G^*$ de rang maximal, et
soit $W' = W_1 \times W_2$ son groupe de Weyl, o\`u $W_1$ et $W_2$
sont des facteurs simples de type classique. Soient
$(\alpha_1,\beta_1)$ et $(\alpha_2,\beta_2)$ deux paires de partitions
param\'etrant des repr\'esentations irr\'eductibles de $W_1$ et $W_2$
respectivement. Alors toute composante irr\'eductible qui intervient
dans la repr\'esentation induite $\Ind_{W'}^W (\alpha_1,\beta_1)
\boxtimes (\alpha_2,\beta_2)$ est param\'etr\'ee par une paire
$(\alpha,\beta)$ v\'erifiant
\[
\alpha \leq \alpha_1 \vee \alpha_2
\qquad\text{et}\qquad
\beta \leq \beta_1 \vee \beta_2.
\]
De plus, il existe un terme dans la repr\'esentation induite pour
lequel ces in\'egalit\'es sont des \'egalit\'es.
\end{lem}

\section{Repr\'esentations et classes de conjugaison de $A(\orb)$}
\label{sect:ab}

Soit $\Nb_G$ (resp.~$\Nb'_G$) l'ensemble des couples $(\orb,\pi)$
(resp.~$(\orb,C)$) o\`u $\orb$ est une classe unipotente et $\pi$ est
une repr\'esentation irr\'eductible (resp.~$C$ est une classe de
conjugaison) de $\Ab(\orb)$. Dans \cite{A} a \'et\'e introduit un
ordre partiel naturel sur $\Nb'_G$. Dans cette section nous traduisons
cet ordre en un ordre partiel sur $\Nb_G$, dans le but de comprendre
la propri\'et\'e d'\^etre bien support\'ee en termes de cette ordre
partiel. En \cite{GGGG}, Geck \`a \'etabli une suite de propri\'et\'es
des faisceaux caract\`eres associ\'es \`a des repr\'esentations
induites v\'erifiant une condition particuli\`ere \`a l'\'egard de
l'ordre partiel usuel sur l'ensemble des classes unipotentes.  La
compr\'ehension de la relation entre l'ordre partiel sur $\Nb_G$ et la
propri\'et\'e d'\^etre bien support\'e devrait permettre d'\'etendre
ces r\'esultats de Geck au cas o\`u sa condition n'est pas satisfaite.

Nous allons \'etudier les groupes $\Ab(\orb)$ comme groupes de
Coxeter, d'apr\`es des id\'ees de \cite{LNU} et de \cite{AS}. Chacun
des groupes $\Ab(\orb)$ est soit un produit de plusieurs exemplaires
de $S_2$ soit l'un des groupes $S_3$, $S_4$, ou $S_5$. Donc de tels
groupes sont tous des groupes de Weyl de type $A$. Pour la discussion
suivante, $H$ d\'esignera un produit quelconque de groupes de Weyl de
type $A$. Supposons choisi un ensemble $\Pi$ de r\'eflexions simples
pour $H$. On peut associer \`a tout sous-ensemble $P$ de $\Pi$ le
sous-groupe $H_P$ de $H$ engendr\'e par les \'el\'ements de
$P$. Puisque tous ces groupes sont des produits de certains groupes
sym\'etriques, il est facile d'\'etablir la proposition suivante.

\begin{prop}\label{prop:bij-irr-p}
Il y a une bijection entre l'ensemble $\Irr(H)$ des
repr\'esentations irr\'eductibles de $H$ et l'ensemble des
sous-ensembles de $\Pi$ \`a conjugaison pr\`es, donn\'ee par
\[
\pi \leftrightsquigarrow P \qquad
\text{si}\qquad
\pi \simeq \epsilon \otimes j_{H_P}^{H} \epsilon,
\]
o\`u $j$ d\'esigne l'induction tronqu\'ee et $\epsilon$ est la
repr\'esentation signe.
\end{prop}

\noindent
(La raison pour laquelle nous tensorisons la repr\'esentation induite
avec la repr\'esentation signe est que plus tard, nous pr\'ef\'ererons
que la repr\'esentation triviale soit la plus grande dans l'ordre
partiel. Si nous ne tensorisions pas ici, la repr\'esentation signe
serait la plus grande.)

Nous d\'ecrivons ensuite une recette pour associer une classe de
conjugaison dans $H$ \`a tout sous-ensemble $P$ de $\Pi$: si
$P = \{s_1, \ldots, s_k\}$, notons $C_P$ la classe de conjugaison de
l'\'el\'ement $s_1 \cdots s_k \in H$. Ce dernier \'el\'ement
d\'epend, bien s\^ur, sur l'ordre dans lequel nous avons \'ecrit les
\'el\'ements de $P$, mais gr\^ace au fait que $H$ est un produit
de certains groupes sym\'etriques, il est facile de d\'emontrer que la
classe $C_P$ n'en d\'epend pas. (On commence par se rappeler que les
classes de conjugaison dans $S_n$ sont param\'etr\'ees par les
partitions de $n$, de sorte que les parties d'une partition donnent
les longueurs des cycles faisant partie d'un \'el\'ement de la classe.
Ensuite on remarque que cette partition peut \^etre calcul\'ee
d'apr\`es une liste non ordonn\'ee de r\'eflexions simples).

Il est \'evident que si $P$ et $Q$ sont deux sous-ensembles de $\Pi$,
alors les classes $C_P$ et $C_Q$ sont \'egales si et seulement si $P$
et $Q$ sont conjugu\'ees.  Donc la recette ci-dessus nous donne
\'egalement une fa\c con d'associer une classe de conjugaison \`a
toute repr\'esentation.

\begin{prop}\label{prop:cl-irr-bij}
L'application $\pi \mapsto P \mapsto C_P$ est une bijection entre
$\Irr(H)$ et l'ensemble $\Cl(H)$ des classes de conjugaison de
$H$.
\end{prop}

Pourtant, cette application n'est pas du tout canonique: elle d\'epend
sur le choix des r\'eflexions simples.  N\'eanmoins, nous pourrons
r\'esoudre ce probl\`eme en utilisant l'ordre partiel sur $\Nb'_G$
pour obtenir un ensemble canonique de r\'eflexions simples.  Rappelons
que l'ensemble $\Cl(H)$ h\'erite d'un ordre partiel de $\Nb'_G$ de
sorte que la classe triviale soit l'unique \'el\'ement minimal.
Appelons une classe $C$ \emph{superminimale} si elle a la
propri\'et\'e que $C > C'$ implique que $C'$ est triviale: les
classes superminimales sont les classes aussi petites que possible
sans \^etre triviales. Le choix des r\'eflexions simples dans
\cite{AS} consiste pr\'ecis\'ement en des repr\'esentants de classes
de conjugaison superminimales.

Le r\'esultat suivant a \'et\'e obtenu de mani\`ere ind\'ependante par
Sommers \cite{S}.

\begin{prop}\label{prop:ab-coxeter}
Il existe un ensemble $\Sigma\subset\Ab(\orb)$ d'involutions, unique
modulo conjugaison, tel que
\begin{enumerate}
\item tout \'el\'ement de $\Sigma$ appartient \`a une classe de
conjugaison superminimale; 
\item 
toute classe de conjugaison superminimale
poss\`ede au moins un repr\'esentant dans $\Sigma$;
\item 
$\Sigma$ constitue un ensemble de r\'eflexions simples pour la
pr\'esentation de $\Ab(\orb)$ comme groupe de Coxeter.
\end{enumerate}
\end{prop}
\begin{proof}[D\'emonstration pour les groupes exceptionnels]
Si $\orb$ est une classe unipotente dans une groupe exceptionnel avec
$A(\orb)$ non trivial, on a toujours que $\Ab(\orb) \simeq S_n$ avec
$2 \leq n \leq 5$. Un coup d'\oe il sur les tables de \cite{A} montre
qu'il y a toujours une unique classe superminimale, dans laquelle se
trouvent toutes les r\'eflexions du groupe. Il est donc \'evident que
l'on peut choisir un ensemble de r\'eflexions simples de sorte que les
conditions ci-dessus soient satisfaites.
\end{proof}

Avant de prouver la proposition pr\'ec\'edente pour les groupes
classiques, nous avons besoin d'une description plus pr\'ecise de
l'ordre partiel sur $\Cl(\Ab(\orb))$. La proposition suivante, qui
fournit une telle description, emploie le langage des \emph{partitions
marqu\'ees}, qui sont des objets combinatoires param\'etrant
l'ensemble $\Nb'_G$, ainsi que la d\'efinition de l'ordre partiel sur
$\Nb'_G$ en termes de la \emph{dualit\'e g\'en\'eralis\'ee} de
Sommers, not\'ee $d_S$ (voir \cite{A}). 

\begin{prop}
Soient $\cmp{\lambda}{\nu}$ et $\cmp{\lambda}{\nu'}$ deux partitions
marqu\'ees qui sont associ\'ees \`a la m\^eme classe unipotente.
Supposons que $\nu$ et $\nu'$ ont tous deux un nombre pair de parties,
en ajoutant un $0$ \`a la fin dans type $C$ si n\'ecessaire. 
\'Ecrivons les parties des partitions marquantes comme 
\begin{align*}
\nu  &= [ a_1 < b_1 < \cdots < a_k < b_k ] \\
\nu' &= [ c_1 < d_1 < \cdots < c_l < d_l ].
\end{align*}
Alors $\cmp{\lambda}{\nu} \geq \cmp{\lambda}{\nu'}$ si et seulement
si, pour tout $i$, $1 \leq i \leq l$, il existe un $j$, $1 \leq j\leq
k$, tel que
\begin{equation}\label{eqn:bloc}
a_j \leq c_i < d_i \leq b_j.
\end{equation}
\end{prop}
\begin{proof}
Si $\nu$ et $\nu'$ v\'erifient la condition ci-dessus, il est facile
de voir que $\cmp{\lambda}{\nu} \geq \cmp{\lambda}{\nu'}$ \`a partir
des formules de \cite{A}, par la m\'ethode de ``division en blocs.''
En particulier, la condition~(\ref{eqn:bloc}) implique qu'une division
en blocs de $\cmp{\lambda}{\nu}$ est toujours une division en blocs de
$\cmp{\lambda}{\nu'}$. La preuve se ram\`ene, alors, au cas o\`u
$\cmp{\lambda}{\nu}$ est un ``bloc de base'': c'est-\`a-dire, $\nu$
n'a que deux parties, qui sont respectivement la plus grande partie
marquable de $\lambda$ et la plus petite partie de $\lambda$ (qui est
suppos\'ee marquable). Dans ce cas, une comparaison directe des
formules de \cite[Proposition 4.12]{A} et celles pour $\dS$ montre
que l'on a toujours $\dS(\cmp{\lambda}{\nu}) \leq
\dS(\cmp{\lambda}{\nu'})$, et donc $\cmp{\lambda}{\nu} \geq
\cmp{\lambda}{\nu'}$.

Pour l'autre implication, remarquons que si la
condition~(\ref{eqn:bloc}) est satisfaite, toute partie de $\nu'$ de
hauteur impaire a une hauteur g\'en\'eralis\'ee impaire dans
$\nu$. D'autre part, toute partie de $\nu$ de hauteur paire a une
hauteur g\'en\'eralis\'ee paire dans $\nu$. Par contre, si la
condition~(\ref{eqn:bloc}) n'est pas satisfaite, alors il existe ou
une partie de $\nu'$ de hauteur impaire dont la hauteur
g\'en\'eralis\'ee dans $\nu$ est paire, ou une partie de $\nu$ de
hauteur paire dont la hauteur g\'en\'eralis\'ee dans $\nu'$ est
impaire. Choisissons une telle partie.  Nous pouvons maintenant
r\'ep\'eter le calcul qui est effectu\'e au cours de la preuve du
Theorem~5.1 de \cite{A}, en faisant jouer le r\^ole de $a$ notre
partie choisie. Ce calcul montre que $\dS(\cmp{\lambda}{\nu}) \not\leq
\dS(\cmp{\lambda}{\nu'})$. Nous en concluons que si~(\ref{eqn:bloc})
n'est pas satisfaite, alors $\cmp{\lambda}{\nu} \not\geq
\cmp{\lambda}{\nu'}$.
\end{proof}

\begin{proof}[D\'emonstration de la Proposition~\ref{prop:ab-coxeter}
pour les groupes classiques] Dans ce cas, le\break groupe $\Ab(\orb)$ est
toujours un produit de plusieurs exemplaires de $S_2$, donc tout
\'el\'ement est une r\'eflexion. Puisque le groupe est ab\'elien,
toute classe de conjugaison ne contient qu'un seul \'el\'ement. Donc
il faut simplement v\'erifier que les \'el\'ements superminimaux
constituent un ensemble de r\'eflexions simples.  Si $a_1 < \cdots <
a_k$ sont les parties marquables d'une partition $\lambda$, la
proposition pr\'ec\'edente implique que les classes superminimales
sont
\[
\cmp{\lambda}{[a_1,a_2]}, \cmp{\lambda}{[a_2,a_3]}, \ldots,
\cmp{\lambda}{[a_{k-1},a_k]}\quad\text{(ainsi que
$\cmp{\lambda}{[a_1]}$ dans type $C$).}
\]
D'apr\`es \cite{LIC}, on sait regarder $A(\orb)$ comme sous-quotient
d'un $S_2$-module libre engendr\'e par $k$ g\'en\'erateurs. On en
d\'eduit et que les \'el\'ements nomm\'es ci-dessus engendrent
$\Ab(\orb)$, et que le nombre d'\'el\'ements est \'egal au rang de
$\Ab(\orb)$ comme $S_2$-module. Par cons\'equent, nous avons
trouv\'e l'ensemble $\Sigma$ tel que d\'ecrit dans la proposition.
\end{proof}

Les Propositions~\ref{prop:cl-irr-bij} et~\ref{prop:ab-coxeter}
ensemble nous fournissent une bijection naturelle $\Nb_G
\leftrightsquigarrow \Nb'_G$. Enfin la relation promise au d\'ebut de
la section est pr\'ecis\'e dans le th\'eor\`eme suivant. \`A partir de
la D\'efinition~\ref{defn:b-s}, ce th\'eor\`eme est presque trivial,
mais peut-\^etre que son \'enonc\'e aide \`a
\'eclaircir ce que signifie \^etre bien support\'e.

\begin{thm}
Les termes d'une repr\'esentation bien support\'ee d'un groupe de Weyl
$W$ appartenant \`a la cellule bilat\`ere la plus haute
correspondent tous \`a des \'el\'ements de $\Nb_G$. De plus, parmi les
termes qui sont associ\'es \`a l'unique classe maximale dans le
support, il y a un unique terme maximal \`a l'\'egard de l'ordre
partiel naturel sur $\Nb_G$.
\end{thm}
\begin{proof}
Le premier \'enonc\'e \'equivaut au Lemme~\ref{lem:ab-defini}. Le
deuxi\`eme \'enonc\'e est cons\'equence du fait que la
repr\'esentation de Springer de l'unique classe maximale $\orb_0$ dans
le support (qui a toujours lieu dans une repr\'esentation bien
support\'ee) correspond \`a la repr\'esentation triviale de
$\Ab(\orb_0)$, qui est toujours maximale dans l'ordre partiel sur
$\Cl(\Ab(\orb_0))$.
\end{proof}

\section{Application aux supports unipotents des faisceaux caract\`eres}
\label{sect:su-fc}

Les faisceaux caract\`eres sur $G$ sont certains faisceaux pervers
$G$-\'equivariants dans la cat\'egorie d\'eriv\'ee des
$\Qlb$-faisceaux constructibles sur $G$, qui ont \'et\'e introduits
par Lusztig dans \cite{LCS}. Nous notons $\Gh$ leur ensemble et
rappelons bri\`evement quelques r\'esultats concernant leur
classification.

Nous supposons que $p$ est bon pour $G$, que le centre $\Centre(G)$ de
$G$ est connexe et que $G/\Centre(G)$ est simple.  Soit $T$ un tore
maximal de $G$ fix\'e et soit $W=\Norm_G(T)/T$ le groupe de Weyl de
$G$ associ\'e. Soient $G^*$ le dual de Langlands de $G$ et $T^*\subset
G^*$ un tore maximal dual de $T$. Il existe une surjection canonique
de $\Gh$ sur l'ensemble des $W$-orbites sur $T^*$,
\cite[Corollaire~11.4]{LCS}. Soient $s\in T^*$ et $(s)$ son orbite
sous $W$. Nous notons $\Gh_s$ l'ensemble des faisceaux caract\`eres
qui appartiennent \`a la fibre au-dessus de $(s)$ de la surjection
cit\'ee. Soit $W_s$ le groupe de Weyl relativement \`a $T^*$ du
centralisateur $G_s^*=\Cent_{G^*}(s)$ de $s$ (ce dernier groupe est
connexe, puisque $\Centre(G)$ l'est).
 
Rappelons maintenant que la correspondance de Springer
g\'en\'eralis\'ee, due \`a Lusztig, est une application qui \'etend la
correspondance de Springer d'origine en une bijection
\[
\nu\colon\bigsqcup\Irr(W_L^G)\,\isom\, \N_G,
\]
l'union \'etant prise sur les classes de $G$-conjugaison de couples
$(L,\iota_0)$, o\`u $L$ est un sous-groupe de Levi de $G$ et $\iota_0$
un \'el\'ement ``cuspidal'' de $\N_L$, et o\`u $W_L^G=\Norm_G(L)/L$
(voir \cite{LIC} et \cite[\S 4.4]{LUS}). En particulier, cette
application attache \`a tout couple $(\orb,\pi) \in \N_G$ un certain
sous-groupe de Levi $L$.

Soit $\orb$ une classe unipotente dans $G$ et soit $\E$ un syst\`eme
local irr\'eductible $G$-\'equivariant sur $\orb$. Nous identifions le
couple $(\orb,\E)$ avec le couple correspondant $\iota =
(\orb,\pi)\in\N_G$ et notons $\IC(\bar\orb,\pi)$ le complexe de
cohomologie d'intersection sur l'adh\'erence de Zariski de $\orb$
associ\'e.  Nous posons $d_\iota=\dim\orb+\dim\Centre(L)$, o\`u $L$
est le sous-groupe de Levi attach\'e \`a $\iota$ par la correspondance
de Springer g\'en\'eralis\'ee.  D'apr\`es \cite[(2.6)(e)]{LCV} et
\cite[6.5]{LIC}, la restriction d'un faisceau caract\`ere $A\in\Gh$
\`a la vari\'et\'e unipotente $G_{\unip}$ de $G$ s'\'ecrit:
\begin{equation} \label{eqn:resunipfc}
A|_{G_{\unip}}=\sum_{\iota\in\N_G}m_{A,\iota}\,A_\iota, \;\;\text{o\`u
$A_\iota=\IC(\bar\orb,\pi)[d_\iota]$,}
\end{equation}
et o\`u les $m_{A,\iota}$ sont certains entiers naturels. Si la
restriction de $A$ \`a $G_{\unip}$ est non nulle, alors il existe au
moins un $\iota$ tel que $m_{A,\iota}\ne 0$. Nous allons d\'ecrire ces
entiers $m_{A,\iota}$: c'est l'objet de la
Proposition~\ref{prop:multiplicites}.
 
\smallskip

Soit $A$ un faisceau caract\`ere sur $G$ dont la restriction \`a
$G_{\unip}$ est non nulle. Il peut \^etre obtenu comme facteur direct
d'un ``induit parabolique'' $\ind_L^G(A_0)$ d'un faisceau caract\`ere
cuspidal $A_0$ sur un sous-groupe de Levi $L$ d'un sous-groupe
parabolique de $G$, tel que la restriction de $A_0$ \`a $L_{\unip}$
est non nulle (ici $\ind_L^G$ d\'esigne l'induction des faisceaux
caract\`eres d\'efinie en \cite[\S 4]{LCS}). Il existe un \'el\'ement
semi-simple $s\in T^*$ tel que $A_0\in\Lh_s$ et les faisceaux
caract\`eres intervenant dans $\ind_L^G(A_0)$ sont param\'etr\'es par
les classes d'isomorphisme de repr\'esentations irr\'eductibles du
pseudo sous-groupe parabolique $W_{L,s}=\Norm_{G^*_s}(L^*)/L^*_s$ du
dual du groupe de Coxeter fini $W_L^G$.  Nous notons $A_{E_1}^s$ le
faisceau caract\`ere d\'efini par la repr\'esentation
$E_1\in\Irr(W_{L,s})$.

\smallskip

Lorsque $G$ est un groupe exceptionnel, $L$ est soit un tore soit le
groupe $G$ lui-m\^eme, et $W_L^G$ est donc soit $W$ soit le groupe
trivial. Soit $G_L$ le groupe de m\^eme type que $G$ de groupe de Weyl
$W_L^G$ et $\N_{G_L}^{\ord}$ l'image de $\Irr(W^G_L)$ dans $\N_{G_L}$
par la correspondance de Springer (non g\'en\'eralis\'ee) pour le
groupe $G_L$. \`A tout $u$-symbole est attach\'e un entier relatif
$d$, appel\'e le \emph{d\'efaut} du $u$-symbole. Dans le cas des
groupes classiques, il r\'esulte de la classification des faisceaux
caract\`eres cuspidaux que $W_L^G$ est de la forme suivante et
l'ensemble $\N_{G_L}^{\ord}$ est en bijection \`a la fois avec $\Irr(W_L^G)$ et
avec le sous-ensemble de l'ensemble des $u$-symboles param\'etrant
$\N_G$ form\'e des $u$-symboles de d\'efaut $d$, o\`u $d$ est d\'efini
comme suit.
\begin{center}
\begin{tabular}{lr@{}lll}
&&\hfil $W_L^G$&&\hfil D\'efaut\\
Type $B$~: & & $B_{n-2t^2-2t}$     & & $d=1+2t$\\
Type $C$~: & & $C_{n-(8t^2\pm2t)}$ & & $d= 1\pm 4t $\\
Type $D$~: & $\bigg\{$ &
\begin{tabular}{@{}l@{}} $D_n$ \\ $B_{n-8t^2}$ \end{tabular} &
\begin{tabular}{@{}l@{}} si $L=T$ \\ sinon \end{tabular} &
\begin{tabular}{@{}l@{}} $d=0$ \\ $d = 4t$ \end{tabular} \\
\end{tabular},
\end{center}
o\`u $t$ est un certain entier positif (nul si et seulement si $L=T$).

Soit $\gamma\colon \N_{G_L}^{\ord}\to\N_G$ le compos\'e
\[
\N_{G_L}^{\ord} \xleftarrow{\simeq} W^G_L \xrightarrow{\nu}\N_G
\]
qui est \'egal \`a l'inclusion canonique
$\N_G^{\ord}\hookrightarrow\N_G$, lorsque $L$ est \'egal \`a
$T$. Lusztig a donn\'e des formules pour $\gamma$ en termes des
$u$-symboles lorsque $G$ est un groupe classique et $L \ne T$ \cite[\S
12.2 et \S 13.2]{LIC}, mais sa description n'est pas correcte dans le
cas o\`u $G$ est de type $C$ et l'image de $\N_{G_L}$ consiste en des
$u$-symboles de d\'efaut n\'egatif. Shoji a expliqu\'e la correction
dans \cite[Remark~5.8]{ShU}: dans ce cas-l\`a, il faut employer une
certaine bijection entre $\Irr(W^G_L)$ et les $u$-symboles de d\'efaut
$1$ diff\'erente que celle d\'ecrite \`a la
Section~\ref{sect:biensupp}. (\emph{Caveat lector}: la formule~(5.4.2)
de \cite{ShU}, qui aurait d\^u \^etre \'egale \`a celle de Lusztig
selon le Remark~5.8, ne l'est pas. Nous remercions F.~L\"ubeck pour
nous avoir indiqu\'e la correction: au cas de d\'efaut positif, il
faut remplacer l'expression ``$B+(2d-1)$'' par ``$B+(2d-2)$'').
Cependant, pour nos propres calculs \`a venir, il sera plus commode de
conserver cette derni\`ere bijection, et faire plut\^ot la
modification au niveau des $u$-symboles. Les formules suivantes pour
$\gamma\binom{\bar\xi_L}{\bar\eta_L}$, dans ce dernier cadre, se
d\'eduisent imm\'ediatement de celles de Shoji. Nous avons \'ecrit
$\bar\xi_L = [a_1 < \cdots < \text{$a_m$ (ou $a_{m+1}$)}]$ et
$\bar\eta_L = [b_1 < \cdots < b_m]$.
\begin{equation}\label{eqn:orb1}
\begin{array}{lc}
\text{Type $B$:} &
\begin{ssymbole}
0 \quad 2 \quad \cdots \quad 4t-2 \quad
 a_1+4t \quad \cdots \quad a_{m+1}+4t\\ b_1 \quad \cdots \quad b_m
\end{ssymbole} \\
\text{Type $C$:} &
\begin{cases}
\hfil\begin{ssymbole}
0 \quad 2 \quad \cdots \quad 8t-2 \quad
 a_1+8t \quad \cdots \quad a_{m+1}+8t\\ b_1 \quad \cdots \quad b_m
\end{ssymbole} &
\rlap{\text{si $d \geq 1$}}\\
\hfil\begin{ssymbole}
b_1 \quad \cdots \quad b_m \\ 1 \quad 3 \quad \cdots \quad
  8t-5 \quad a_1 + 8t-3 \quad \cdots \quad a_{m+1} + 8t-3
\end{ssymbole} &
\rlap{\text{si $d \leq -1$}}
\end{cases} \\
\text{Type $D$:} &
\begin{ssymbole}
0 \quad 2 \quad \cdots \quad 8t-4 \quad
 a_1+8t-2 \quad \cdots \quad a_{m+1}+8t-2\\ b_1 \quad \cdots \quad b_m
\end{ssymbole}
\end{array}
\end{equation}

\smallskip

Soient $L_{\ad}$ et $L_{\der}$ respectivement le groupe adjoint et le
sous-groupe d\'eriv\'e de $L$ et soit $\pr \colon L \rightarrow
L_{\ad}$ la projection canonique. D'apr\`es \cite[(17.10)]{LCS}, il
existe un faisceau caract\`ere cuspidal $\bA_0$ sur $L_{\ad}$ tel que
$A_0=\pr^{*}(\bA_0) \otimes \CL$, o\`u $\CL$ est un syst\`eme local de
Kummer sur $L$, qui est l'image r\'eciproque d'un syst\`eme local sur
$L/L_{\der}$ sous l'application canonique $L\rightarrow L/L_{\der}$.
Soit $L_{\der}^* \hookrightarrow L^*$ le plongement correspondant
entre les groupes duaux.  Si $\bA_0$ appartient \`a
$(\widehat{L_{\ad}})_{\bs}$, avec $\bs\in T^* \cap L_{\der}^*$, et
$\CL$ correspond \`a l'\'el\'ement central $z$ de $L^*$, alors
$A_0\in\Lh_s$, o\`u $s=\bs z$. Remarquons que $L^*_s=L^*_{\bs}$.

Les faisceaux caract\`eres sur $G$ qui interviennent dans
$\ind_L^G(\pr^{*}(\bA_0))$ sont eux param\'etr\'es par les classes
d'isomorphisme de repr\'esentations irr\'eductibles du groupe
$W_L^G=\Norm_{G^*}(L^*)/L^*\simeq W_{L,\bs}$.  Nous noterons
$A^{\bs}_{E'}$ le faisceau caract\`ere d\'efini par
$E'\in\Irr(W_L^G)$. La restriction de $A^{\bs}_{E'}$ \`a $G_{\unip}$
est \`a support sur $\overline\orb$ et \'egale \`a $A_\iota$ sur sur
$\overline\orb$, si $\iota=(\orb,\pi)=\nu(E')$.  Les faisceaux
caract\`eres $A_0$ et $\pr^{*}(\bA_0)$ ont m\^eme restriction \`a
$L_{\unip}$ et l'on a
\begin{equation} \label{eqn:suppfc}
A_{E_1}^s=\bigoplus_{E'}\, (E':\Ind_{W_{L,s}}^{W_L^G}(E_1))\,A_{E'}^{\bs}
=\bigoplus_{E'}\, (E':\Ind_{W_{L,s}}^{W_L^G}(E_1))\,A_{\nu(E')}
\;\;\text{ sur $G_{\unip}$,}
\end{equation}
o\`u $E'$ parcourt les repr\'esentations irr\'eductibles de $W_L^G$,
\`a isomorphisme pr\`es, et $(E':\Ind_{W_{L,s}}^{W_L^G}(E_1))$
d\'esigne la multiplicit\'e de $E'$ dans la repr\'esentation induite
$E=\Ind_{W_{L,s}}^{W_L}(E_1)$, \cite[(2.6)]{LCS}.

\smallskip
La proposition suivante se d\'eduit imm\'ediatement
de~(\ref{eqn:resunipfc}) et~(\ref{eqn:suppfc}).

\begin{prop} \label{prop:multiplicites}
Soit $A=A_{E_1}^s$ un faisceau caract\`ere sur $G$, avec
$E_1\in\Irr(W_{L,s})$. Alors
\[
m_{A,\iota}=\begin{cases}
(E':\Ind_{W_{L,s}}^{W_L^G}(E_1))&\text{si $\iota=\nu(E')$, avec
$E'\in\Irr(W_L^G)$,}\cr
0&\text{si $\iota\notin\nu(\Irr(W_L^G))$.}\end{cases}
\]
\end{prop}

\smallskip

\begin{thm} \label{thm:supportfc}
Soit $A$ un faisceau caract\`ere non identiquement nul sur la
vari\'et\'e unipotente. Il existe alors une classe unipotente $\orb_A$
sur laquelle la restriction de $A$ n'est pas nulle et dont
l'adh\'erence de Zariski contient toute classe unipotente sur laquelle
cette restriction n'est pas nulle.
\end{thm}
\begin{proof}
La relation~(\ref{eqn:suppfc}) permet d'associer \`a $A$ un certain
sous-ensemble $P \subset \N_{G_L}^{\ord}$, qui contient les
\'el\'ements correspondant \`a des repr\'esentations $E'$ de $W_L^G$
telles que $(E' : \Ind_{W_{L,s}}^{W_L^G} (E_1)) \neq 0$. Dans le cas
o\`u $L = T$, la proposition~\ref{prop:ind-supp} dit que tous les
\'el\'ements de ce sous-ensemble sont associ\'es \`a des classes
unipotentes contenues dans l'adh\'erence de Zariski de $\orb$, et donc
le th\'eor\`eme est d\'emontr\'e. Si $L \neq T$, nous voudrions que
l'image $\gamma(P) \subset \N_G$ ait cette derni\`ere
propri\'et\'e. Nous traiterons les cas de $G$ classique et de $G$
exceptionnel s\'epar\'ement.

Pour $G$ exceptionnel, la seule possibilit\'e est que $L = G$, comme
nous avons d\'ej\`a remarqu\'e. Dans ce cas l'ensemble
$\N_{G_L}^{\ord}$ ne contient qu'un seul \'el\'ement, et donc il est
\'evident que $\gamma(P)$ a la propri\'et\'e cherch\'ee.

Pour $G$ classique, nous pouvons appliquer le
Lemme~\ref{lem:ind-part}. La conclusion de ce lemme-l\`a, traduite
dans le langage des $u$-symboles, dit qu'il existe un unique
\'el\'ement de $P$ dont le $u$-symbole $\binom{\bar\xi_0}{\bar\eta_0}$
a la propri\'et\'e que
\begin{equation}\label{eqn:rel}
\bar\xi \leq \bar\xi_0
\qquad\text{et}\qquad
\bar\eta \leq \bar\eta_0
\end{equation}
pour tout autre $u$-symbole $\binom{\bar\xi}{\bar\eta}$ qui est membre
de $P$. La description ci-dessus de $\gamma$ montre que $\gamma$
conserve la relation~(\ref{eqn:rel}), et donc il est possible
d'appliquer le Lemme~\ref{lem:ordre} aux \'el\'ements de $\gamma(P)$.
Le th\'eor\`eme en d\'ecoule.
\end{proof}

\begin{defn} \label{defn:supportfc}
La classe $\orb_A$ est appel\'ee le \emph{support unipotent} de $A$.
\end{defn}

Lusztig a d\'efini une surjection canonique de $\Gh_s$ sur l'ensemble
des cellules bilat\`eres de $W_s$ en \cite[Corollary~16.7]{LCS}. Si
$\bc$ est une cellule bilat\`ere donn\'ee de $W_s$, nous notons
$\Gh_{s,\bc}$ l'ensemble des faisceaux-caract\`eres dans la fibre
au-dessus de $\bc$ de cette surjection.  Nous obtenons ainsi une
partition de $\Gh_s$:
\[
\Gh_s=\bigsqcup_{\bc}\Gh_{s,\bc},
\]
o\`u $\bc$ parcourt les cellules bilat\`eres de $W_s$.  Nous associons
une classe unipotente $\orb_{s,\bc}$ \`a l'ensemble $\Gh_{s,\bc}$ de
la mani\`ere suivante: soit $E_1$ la repr\'esentation sp\'eciale de
$W_s$, d'apr\`es le Th\'eor\`eme~\ref{thm:b-s}, la repr\'esentation
induite $E=\Ind_{W_s}^W(E_1)$ est bien support\'ee et $\orb_{s,\bc}$
est d\'efinie comme \'etant l'unique classe maximale dans le support
de $E$.  La preuve du Th\'eor\`eme~\ref{thm:b-s} montre que la classe
$\orb_{s,\bc}$ est la classe induite de la classe associ\'ee \`a $E_1$
par la correspondance de Springer, il s'ensuit que $\orb_{s,\bc}$ est
\'egale \`a la classe unipotente attach\'e \`a $\Gh_{s,\bc}$ par Lusztig
en \cite[(13.3)]{Lbook} ou \cite[\S 10.5]{LUS}.

\begin{prop} \label{prop:incl_supp}
Soit $A$ un faisceau caract\`ere appartenant \`a $\Gh_{s,\bc}$. Alors
la classe $\orb_A$ est contenue dans l'adh\'erence de Zariski de
$\orb_{s,\bc}$.
\end{prop}

\begin{remark}
Il est possible de montrer que $\dim \orb_A \leq \dim \orb_{s,\bc}$
sans faire les longs calculs qui suivent. En effet cette in\'egalit\'e
r\'esulte imm\'ediatement du Th\'eor\`eme~\ref{thm:supportfc},
combin\'e avec \cite[Theorem~1.1.(a)]{Au}.
\end{remark}

\begin{proof}
Lorsque $L=T$, les classes $\orb_A$ et $\orb_{s,\bc}$ sont
\'egales. Nous supposons d\'esormais $L$ diff\'erent de $T$.  La table
ci-dessous r\'esume alors les possibilit\'es pour $W^L_T$ et
$W^L_{T,s}$. De plus, chaque groupe $W^L_{T,s}$ contient une unique
cellule bilat\`ere cuspidale, qui sera not\'ee $\bc_0$. La derni\`ere
colonne de la table donne l'unique repr\'esentation sp\'eciale
$E_{\bc_0}$ dans $\bc_0$ \cite[\S 8.1]{Lbook}, param\'etr\'ee par des
paires de partition. (Ici $\alpha_t$ d\'esigne la partition $[1 < 2 <
\cdots < t]$ de $t(t+1)/2$).
\[
\begin{array}{cr@{}cccl}
    && W^L_T & W^L_{T,s} & E_{\bc_0} \\
\text{Type $B$:} && B_{2t^2+2t} & C_{t^2+t}\times C_{t^2+t} &
(\alpha_t,\alpha_t) \boxtimes (\alpha_t,\alpha_t) \\
\text{Type $C$:} & \bigg\{ &
\begin{array}{@{}c@{}}
C_{8t^2+2t} \\ C_{8t^2-2t}
\end{array}
&
\begin{array}{@{}c@{}}
D_{4t^2}\times B_{4t^2+2t} \\ D_{4t^2}\times B_{4t^2-2t}
\end{array}
&
\begin{array}{@{}c@{}}
(\alpha_{2t},\alpha_{2t-1}) \boxtimes (\alpha_{2t},\alpha_{2t}) \\
(\alpha_{2t},\alpha_{2t-1}) \boxtimes (\alpha_{2t-1},\alpha_{2t-1}) 
\end{array}
&
\begin{array}{@{}l@{}}
\text{si $d \geq 1$}  \\ \text{si $d \leq -1$}
\end{array}
\\
\text{Type $D$:} && D_{8t^2} & D_{4t^2} \times D_{4t^2} &
(\alpha_{2t},\alpha_{2t-1}) \boxtimes (\alpha_{2t},\alpha_{2t-1})
\end{array}
\]

Soit $(\orb_0,\pi_0)$ l'image par $\gamma$ de la repr\'esentation
triviale de $W^G_L$. La preuve de la proposition proc\`ede en deux
\'etapes: nous \'etablissons d'abord que $\orb_0 = \orb_{s,\bc}$, et
puis nous d\'emontrons que pour tout $A$, la classe $\orb_A$ est
contenue dans l'adh\'erence de Zariski de $\orb_0$.

Quel que soit le type de $W^G_L$, sa repr\'esentation triviale est
param\'etr\'ee par le couple $([k],\varnothing)$, o\`u $k$ est son
rang, et donc le $u$-symbole correspondant est
$\binom{k}{\varnothing}$. D'apr\`es~(\ref{eqn:orb1}), l'image par
$\gamma$ de ce $u$-symbole, \ie le $u$-symbole de la paire
$(\orb_0,\pi_0)$, est donc
\begin{equation}\label{eqn:orb0}
\begin{array}{lc}
\text{Type $B$:} &
\begin{ssymbole}
0 \quad 2 \quad \cdots \quad 4t-2 \quad (n-2t^2-2t)+4t \\
\varnothing
\end{ssymbole} \\
\text{Type $C$:} &
\begin{cases}
\hfil\begin{ssymbole}
0 \quad 2 \quad \cdots \quad 8t-2 \quad (n-8t^2-2t)+8t \\
\varnothing
\end{ssymbole} &
\rlap{\text{si $d \geq 1$}} \\
\hfil\begin{ssymbole}
\varnothing \\
1 \quad 3 \quad \cdots \quad 8t-5 \quad (n-8t^2+2t)+8t-3
\end{ssymbole} &
\rlap{\text{si $d \leq -1$}}
\end{cases} \\
\text{Type $D$:} &
\begin{ssymbole}
0 \quad 2 \quad \cdots \quad 8t-4 \quad (n-8t^2)+8t-2 \\
\varnothing
\end{ssymbole}
\end{array}
\end{equation}

Ensuite nous calculons le $u$-symbole de la paire $(\orb_{s,\bc},
1)$. Si aucun facteur de type $D$ n'a lieu dans $W^L_{T,s}$, alors il
faut simplement appliquer le Lemme~\ref{lem:ind-part}: la
repr\'esentation qu'il fournit, param\'etr\'ee par les partitions jointes
des partitions d'origine, est bien celle obtenue par induction
tronqu\'ee. D'autre part, la seule repr\'esentation qu'il faut traiter
dans le cas o\`u il y a un facteur de type $D$ est $(\alpha_{2t},
\alpha_{2t-1})$. Les \'egalit\'es suivantes d\'ecoulent de la
description des groupes de Weyl de type $D$ et les formules pour la
$a$-valeur en \cite{LIR}:
\begin{align*}
\Ind_{D_{4t^2}}^{B_{4t^2}} (\alpha_{2t},\alpha_{2t-1}) &=
(\alpha_{2t},\alpha_{2t-1}) \oplus (\alpha_{2t-1},\alpha_{2t}) \\
j_{D_{4t^2}}^{B_{4t^2}} (\alpha_{2t},\alpha_{2t-1}) &=
(\alpha_{2t},\alpha_{2t-1})
\end{align*}
Gr\^ace \`a la transitivit\'e de l'induction tronqu\'ee, ce fait
permet de finir le calcul en utilisant le Lemme~\ref{lem:ind-part}. La
table suivante donne les r\'esultats de ces calculs. (Ici $2\alpha_t$
d\'esigne la partition jointe $\alpha_t \vee \alpha_t$, et $k$ est
toujours le rang de $W^G_L$).
\[
\begin{array}{lr@{}ccl}
&& j_{W^L_{T,s}}^{W^L_T} E_{\bc_0} & j_{W^L_{T,s}}^{W} E_{\bc_0} \\
\text{Type $B$:} && (2\alpha_t,2\alpha_t) & (
  2\alpha_t \vee [k], 2\alpha_t) \\
\text{Type $C$:} & \bigg\{ & 
\begin{array}{@{}c@{}}
(2\alpha_{2t},\alpha_{2t}\vee\alpha_{2t-1}) \\
(\alpha_{2t}\vee\alpha_{2t-1}, 2\alpha_{2t-1})
\end{array}
&
\begin{array}{@{}c@{}}
(2\alpha_{2t} \vee [k],\alpha_{2t}\vee\alpha_{2t-1}) \\
(\alpha_{2t}\vee\alpha_{2t-1} \vee [k], 2\alpha_{2t-1})
\end{array}
&
\begin{array}{@{}l@{}}
\text{si $d \geq 1$} \\
\text{si $d \leq -1$} \\
\end{array}
\\
\text{Type $D$:} && (2\alpha_{2t},2\alpha_{2t-1}) & 
  (2\alpha_{2t} \vee [k],2\alpha_{2t-1}) \\
\end{array}
\]
Les $u$-symboles correspondants sont alors:
\begin{equation}\label{eqn:orbsc}
\begin{array}{lc}
\text{Type $B$:} &
\begin{ssymbole}
0 &\quad& 4 &\quad& \cdots &\quad& 4t-4 &\quad& (n-2t^2+2t)+4t \\
& 2 && 6 && \cdots && 4t-2
\end{ssymbole} \\
\text{Type $C$:} &
\begin{cases}
\hfil\begin{ssymbole}
0 &\quad& 4 &\quad& \cdots &\quad& 8t-4 &\quad& (n-8t^2-2t)+8t \\
& 2 && 6 && \cdots && 8t-2
\end{ssymbole} &
\rlap{\text{si $d \geq 1$}} \\
\hfil\begin{ssymbole}
1 &\quad& 5 &\quad& \cdots &\quad& 8t-7 &\quad& (n-2t^2+2t)+8t-3 \\
& 3 && 7 && \cdots && 8t-5
\end{ssymbole} &
\rlap{\text{si $d \leq -1$}}
\end{cases} \\
\text{Type $D$:} &
\begin{ssymbole}
2 &\quad& 6 &\quad& \cdots &\quad& 8t-6 &\quad& (n-8t^2)+(8t-2) \\
0 && 4 && \cdots && 8t-8 && 8t-4
\end{ssymbole}
\end{array}
\end{equation}
Il est \'evident que tout $u$-symbole qui parait dans~(\ref{eqn:orb0})
appartient \`a m\^eme classe de similitude que le $u$-symbole
correspondant dans~(\ref{eqn:orbsc}). Nous en concluons que $\orb_0 =
\orb_{s,\bc}$.

La restriction d'un faisceau caract\`ere $A$ intervenant dans
$\ind_L^G A_0$ \`a une classe unipotente $\orb$ est non nulle s'il
existe une paire $(\orb,\pi)$ qui fait partie de l'image par $\gamma$
de l'ensemble $P \subset \N_{G_L}^{\ord}$ associ\'e \`a $A$ comme dans
la preuve du Th\'eor\`eme~\ref{thm:supportfc}. Nous allons montrer ici que
si $(\orb,\pi)$ est l'image d'un \'el\'ement quelconque de
$\N_{G_L}^{\ord}$, alors $\orb \subset \overline{\orb_0}$. Cela aura pour
cons\'equence que $\orb_A$ est contenue dans l'adh\'erence de Zariski
de $\orb_{s,\bc}$.

Soit $\binom{\bar\xi_L}{\bar\eta_L}$ le $u$-symbole d'un \'el\'ement
de $\N_{G_L}^{\ord}$ dont l'image par $\gamma$ est $(\orb,\pi)$.  Le
$u$-symbole de $(\orb,\pi)$, qui sera not\'ee
$\binom{\bar\xi}{\bar\eta}$, et donn\'e dans~(\ref{eqn:orb1}).
Pour faciliter la comparaison de $\orb$ et $\orb_0$, augmentons les
longueurs des $u$-symboles de~(\ref{eqn:orb0}) jusqu'\`a ce qu'elles
soient \'egales \`a celles de~(\ref{eqn:orb1}):
\[
\begin{array}{lc}
\text{Type $B$:} &
\begin{ssymbole}
0 \quad 2 \quad \cdots \quad 4t-2+2m \quad (n-2t^2-2t)+4t+2m \\
0 \quad 2 \quad \cdots \quad 2m-2
\end{ssymbole} \\
\text{Type $C$:} &
\begin{cases}
\hfil\begin{ssymbole}
0 \quad 2 \quad \cdots \quad 8t-2+2m \quad (n-8t^2-2t)+8t+2m \\
1 \quad 3 \quad \cdots \quad 2m-1
\end{ssymbole} &
\rlap{\text{si $d \geq 1$}} \\
\hfil\begin{ssymbole}
0 \quad 2 \quad \cdots \quad 2m-2 \\
1 \quad 3 \quad \cdots \quad 8t-5+2m \quad (n-8t^2+2t)+8t-3+2m
\end{ssymbole} &
\rlap{\text{si $d \geq 1$}}
\end{cases} \\
\text{Type $D$:} &
\begin{ssymbole}
0 \quad 2 \quad \cdots \quad 8t-4+2m \quad (n-8t^2)+8t-2+2m \\
0 \quad 2 \quad \cdots \quad 2m-2
\end{ssymbole}
\end{array}
\]
Notons $\binom{\bar\xi_0}{\bar\eta_0}$ ce dernier $u$-symbole. La
d\'efinition des $u$-symboles impose que le $i$-\`eme entier dans
chaque ligne d'un $u$-symbole soit au minimum $2(i-1)$ ($2i-1$ pour la
deuxi\`eme ligne d'un $u$-symbole de type $C$). Autrement dit, nous
avons \'ecrit les deux d'une telle fa\c con que tout entier dans
$\binom{\bar\xi}{\bar\eta}$, \`a l'exception \'eventuelle du plus
grand entier sur la premi\`ere ligne, soit plus grand que l'entier
correspondant dans $\binom{\bar\xi_0}{\bar\eta_0}$. Soient $\mu =
\bar\xi \cup \bar\eta$ et $\mu_0 = \bar\xi_0 \cup
\bar\eta_0$. L'observation pr\'ec\'edente, traduite dans le langage
de partitions, dit qu'il existe une correspondance entre les parties
de $\mu$ et celles de $\mu_0$ telle que toute partie de $\mu_0$, sauf
peut-\^etre la plus grande, est plus petite que la partie
correspondante de $\mu$. Cela implique que la somme des $i$ parties
les plus petites de $\mu$ est plus grande que la somme correspondante
de certaines parties de $\mu_0$, laquelle est ensuite plus grande que
la somme des $i$ parties les plus petites de $\mu_0$. Puisque la somme
de toutes les parties de $\mu$ est \'egale \`a celle de $\mu_0$, nous
pouvons inverser cette in\'egalit\'e pour conclure que
$\sigma_i(\mu_0) \geq \sigma_i(\mu)$; \ie $\mu_0 \geq \mu$. Le
Lemme~\ref{lem:ordre} implique alors que $\orb \subset
\overline{\orb_0}$.
\end{proof}

Pour tout $A\in\Gh_{s,\bc}$, nous posons
\begin{equation} \label{eqn:supportfc}
\N_G(A)=\left\{(\orb,\pi)=\iota\in \N_G\,:\,\orb=\orb_{s,\bc}\;
\text{ et }\;m_{A,\iota}\ne 0\right\}.
\end{equation}
Cet ensemble appara\^\i t d\'ej\`a dans \cite[Theorem~4.5]{GGGG}.

\begin{cor} \label{cor:supportfc}
Soit $A\in\Gh_{s,\bc}$.  L'ensemble $\N_G(A)$ est non vide si et
seulement si $\orb_A$ est \'egale \`a $\orb_{s,\bc}$.
\end{cor}
\begin{proof}
Le fait que $\N_G(A)$ soit non vide est \'equivalent au fait que la
restriction de $A$ \`a $\orb_{s,\bc}$ soit non nulle. Le
Th\'eor\`eme~\ref{thm:supportfc} montre alors que la classe
$\orb_{s,\bc}$ est contenue dans l'adh\'erence de Zariski de
$\orb_A$. La Proposition~\ref{prop:incl_supp} implique alors que les
classes $\orb_{s,\bc}$ et $\orb$ sont \'egales.
\end{proof}

\smallskip
\begin{remark}
Comme le montre l'exemple suivant, le support unipotent d'un faisceau
caract\`ere $A\in\Gh_{s,\bc}$, non identiquement nul sur la
vari\'et\'e unipotente, peut \^etre diff\'erent de $\orb_{s,\bc}$.
Soit $G$ un groupe de type $F_4$, et soit $s \in G^*$ un \'el\'ement
semi-simple tel que $G^*_s$ soit de type $C_3 \times A_1$. Soit $E_1$
la repr\'esentation de $W_s$ param\'etr\'ee par $([1<2],\varnothing)
\boxtimes [2]$. Cette repr\'esentation n'est pas sp\'eciale. La
repr\'esentation sp\'eciale dans sa cellule bilat\`ere est la
repr\'esentation $([2],[1]) \boxtimes [2]$, dont l'induction
tronqu\'ee est la repr\'esentation de Springer de la classe sp\'eciale
$F_4(a_1)$. Nous avons donc identifi\'e la classe $\orb_{s,\bc}$.
Pourtant, l'image par la correspondance de Springer des
composantes de la repr\'esentation induite $E$ de $E_1$ consiste en
les quatre paires suivantes:
\[
(F_4(a_2),1),\quad (F_4(a_2),\epsilon), \quad
(B_2,1), \quad\text{et}\quad (B_2,\epsilon).
\]
(Ici $\epsilon$ d\'esigne l'unique repr\'esentation non triviale de
$A(\orb)$). En particulier, la classe $\orb_{s,\bc} = F_4(a_1)$ ne
fait pas partie du support de $E$. Le faisceau caract\`ere $A^s_{E_1}$
est donc de restriction nulle sur $\orb_{s,\bc}$.

Remarquons de plus que la classe $F_4(a_2)$ est sp\'eciale, mais que
la repr\'esentation de $W$ qui correspond au syst\`eme local non trivial
appartient \`a la cellule bilat\`ere qui correspond \`a la classe
sp\'eciale $F_4(a_1)$, c'est-\`a-dire, \`a une cellule plus haute que
celle \`a laquelle appartient la repr\'esentation de Springer de
l'unique classe maximale. Donc la repr\'esentation induite de $E_1$
n'est pas bien support\'ee.
\end{remark}

\begin{remark}
Soit $A$ le faisceau caract\`ere intervenant dans $\ind_L^G A_0$ qui
correspond \`a la repr\'esentation triviale de $W^G_{L,s}$.  Les
calculs effectu\'es ci-dessus montrent que $\orb_0$ est l'unique plus
grande classe sur laquelle la restriction de $A$ n'est pas nulle. Ce
faisceau caract\`ere est alors un \'el\'ement de $\Gh_{s,\bc}$ dont le
support unipotent est $\orb_{s,\bc}$. En particulier, ceci montre (en
l'appliquant \`a $L=G$, o\`u le groupe $W^G_{L,s}$ est alors r\'eduit
\`a $\{1\}$) que le support unipotent d'un faisceau caract\`ere
cuspidal appartenant \`a $\Gh_{s,\bc}$ est \'egal \`a $\orb_{s,\bc}$.
\end{remark}

\section{Application aux supports unipotents des caract\`eres}
\label{sect:su-car}

\subsection{Caract\`eres et caract\`eres fant\^omes}

Nous supposons $T$ rationnel et contenu dans un sous-groupe de
Borel rationnel $B$. Le groupe dual $G^*$ h\'erite d'une
structure $\F_q$-rationnelle et nous notons encore $F$ l'endomorphisme
de Frobenius associ\'e \`a celle-ci.  Nous pouvons supposer que les
groupes $B^*$ et $T^*$ sont $F$-stables, \cite{DL}.  
Nous supposons que $F(s)=s$.
Lusztig a associ\'e \`a toute famille $\Fa$ de $W_s$ un groupe fini $\G_\Fa$, 
et d\'ecrit une injection, not\'ee $E\mapsto x_E$,
de la famille $\Fa$ dans l'ensemble fini $\M(\G_\Fa)$ des classes de
$\G_\Fa$-conjugaison de paires $(x,\rho)$, o\`u $x$ est un
\'el\'ement de
$\G_\Fa$ et $\rho$ est une repr\'esentation irr\'eductible sur $\Qlb$
(d\'efinie \`a isomorphisme pr\`es) du centralisateur de $x$ dans
$\G_\Fa$. Chaque famille $\Fa$ contient donc une unique
repr\'esentation sp\'eciale: celle-ci correspond \`a l'\'el\'ement
$(1,1)$ de $\M(\G_\Fa)$. 
L'ensemble $\Gh_{s,\bc}$ est en bijection avec $\M(\G_\Fa)$,
o\`u $\Fa$ est la famille des repr\'esentations irr\'eductibles de $W_s$
correspondant \`a la cellule $\bc$. Remarquons que
$\G_\Fa = \Ab(\orb)$ o\`u $\orb$ est la classe sp\'eciale qui
correspond \`a la famille $\Fa$.
 
D'autre part, Lusztig a d\'efini un
accouplement sur $\M(\G_\Fa)$ par la formule \cite[(4.14.3)]{Lbook}:
\[\left\{(x,\rho),(y,\tau)\right\}=\sum_{\substack{g\in\G_\Fa\\
x\cdot gyg^{-1}=gyg^{-1}\cdot
x}}\frac{\Tr(gxg^{-1},\tau)\cdot\Tr(gyg^{-1},\rho)}
{|\Cent_{\G_\Fa}(x)|\cdot|\Cent_{\G_\Fa}(y)|}.\]
Nous notons $X(W_s)$ l'union disjointe, sur toutes les familles $\Fa$, des
ensembles $\M(\G_\Fa)$ et nous \'etendons $\{\;,\;\}$ \`a $X(W_s)$ en
posant $\{x,x'\}=0$ si $x$ et $x'$ sont dans des familles distinctes.

\smallskip

Il existe une partition de l'ensemble des caract\`eres irr\'eductibles de
$G^F$ en sous-ensembles $\E(G^F)_s$, appel\'es \emph{s\'eries de
Lusztig}, param\'etr\'es par les classes de conjugaison
semi-simples $F$-stables dans $G^*$. 
Pour simplifier l'exposition, nous supposons dans cette sous-section que
le groupe $G$ a un centre connexe et est
d\'eploy\'e sur $\F_q$. L'ensemble $\E(G^F)_s$ est 
alors en bijection avec $X(W_s)$. Nous notons $\rho_x$ le caract\`ere de $G^F$
param\'etr\'e par l'\'el\'ement $x$ de $X(W_s)$. Nous d\'efinissons
formellement une fonction centrale $R_x$, appel\'ee \emph{caract\`ere
fant\^ome} \cite[(4.24.1) et p.~347]{Lbook}, par
\[R_x=\sum_{y\in X(W_s)}\{x,y\}\,\Delta(y)\,\rho_y,\]
o\`u $\Delta\colon\M(\G_\Fa)\to\{\pm 1\}$ est la fonction d\'efinie en
\cite[(4.14)]{Lbook} (la fonction $\Delta$ ne prend la valeur $-1$ que
dans quelques cas o\`u $W_s$ est de type $E_7$ ou $E_8$). 
La matrice $(\{x,y\})$ est unitaire et les $R_x$, pour $x\in X(W_s)$,
forment une base orthonormale du sous-espace de l'espace des fonctions
centrales sur $G^F$ engendr\'e par les caract\`eres appartenant \`a
$\E_s$.

Soit $\bc$ une cellule bilat\`ere dans $W_s$. La paire $(s,\bc)$ d\'efinit
un sous-ensemble $\E(G^F)_{s,\bc}$ de $\E(G^F)_s$, \cite[(8.4.4) et
(6.17)]{Lbook}. De plus, puisque $F$ agit trivialement sur $W_s$,
il existe,
d'apr\`es \cite[Main Theorem 4.23]{Lbook}, aussi une bijection 
de $\E(G^F)_{s,\bc}$ sur $\M(\G_\Fa)$, o\`u $\Fa$
est la famille de $W_s$ correspondant \`a $\bc$.

\smallskip
Un faisceau pervers $A$ sur $G$ est dit \emph{$F$-stable} si $F^*A$ et
$A$ sont isomorphes.
On associe \`a tout couple $(A,\varphi)$,
form\'e d'un faisceau pervers $F$-stable $G$-\'equivariant $A$ sur $G$ et d'un
isomorphisme $\varphi\colon F^*A\isom A$, une fonction centrale
$\chi_{A,\varphi}$ sur $G^F$ (\`a valeur dans $\Qlb$),
appel\'ee la \emph{fonction caract\'eristique} de $A$ associ\'ee \`a
$\varphi$, d\'efinie
par
\[\chi_{A,\varphi}(g):=\sum_i(-1)^i\Tr(\varphi,\CH_g^i(A)),\;\;\text{ $g\in
G^F$,}\]
o\`u $\CH_g^i(A)$ d\'esigne la fibre en $g$ du $i$-\`eme faisceau de
cohomologie $\CH^i(A)$ de $A$ et o\`u l'on note encore $\varphi$
l'application lin\'eaire que $\varphi$ induit sur $\CH^i(A)$.
Pour les faisceaux caract\`eres de restriction non nulle \`a la
vari\'et\'e unipotente il existe un choix canonique pour $\varphi$,
indiqu\'e en \cite[(3.2)]{LCV}. Nous noterons $\chi_A$ la fonction
caract\'eristique correspondante \`a ce choix.

Le lien entre les caract\`eres fant\^omes de $G^F$ et les fonctions
caract\'eristiques de faisceaux caract\`eres $F$-stables sur $G$ est fourni 
par une conjecture de Lusztig. La forme ci-dessous de cette conjecture a
\'et\'e prouv\'ee par Shoji dans \cite{Sh}: pour tout $x\in X(W_s)$, on a 
$R_x=\zeta_{A_x}\chi_{A_x}$, o\`u $A_x\in\Gh_{s,\bc}$ est param\'etr\'e par
$x$ et $\zeta_{A_x}$ est le nombre alg\'ebrique de module un associ\'e
\`a $A_x$ par \cite[Theorem~13.10.(b)]{LCS}. 

Il en r\'esulte, en particulier, que tout caract\`ere irr\'eductible
$\rho$ de $G^F$ qui appartient \`a $\E(G^F)_{s,\bc}$ est combinaison
lin\'eaire de fonctions caract\'eristiques $\chi_{A,\varphi}$, pour des
faisceaux caract\`eres $A\in\Gh_{s,\bc}$. 
Le th\'eor\`eme suivant est alors une cons\'equence \'evidente de la
Proposition~\ref{prop:incl_supp}.

\begin{thm}
Soit $\rho$ un caract\`ere irr\'eductible de $G^F$ appartenant \`a
$\E(G^F)_{s,\bc}$.
Toute classe unipotente rationnelle sur laquelle la restriction de
$\rho$ est non identiquement
nulle est contenue dans l'adh\'erence de Zariski de $\orb_{s,\bc}$.
\end{thm}

Nous allons voir ci-dessous que la restriction de
$\rho$ \`a la classe $\orb_{s,\bc}$ est non identiquement nulle.
Nous d\'efinissons la \emph{valeur moyenne} $\VM(f,\orb)$ et la 
\emph{valeur moyenne pond\'er\'ee} $\VMP(f,\orb)$, 
sur les points $\mathbb{F}_q$-rationnels d'une classe unipotente
$F$-stable $\orb$, d'une fonction centrale $f$ sur le groupe de
Chevalley fini $G(\mathbb{F}_q)$, de la mani\`ere suivante: notons
$u_1$,$u_2$, $\ldots$, $u_r$ des repr\'esentants dans
$\orb(\mathbb{F}_q)$ des classes de $G(\mathbb{F}_q)$-conjugaison
contenues dans $\orb(\mathbb{F}_q)$, alors
$$
\VM(f,\orb):= \sum_{i=1}^r[G^F:\Cent_{G}(u_i)^F]\,f(u_i)
=\sum_{u\in \orb(\mathbb{F}_q)}f(u),
$$ 
$$
\VMP(f,\orb):=\sum_{i=1}^r [A(u_i):A(u_i)^F]\,f(u_i),
$$ 
o\`u $A(u_i)^F$ d\'esigne le groupe des points de $A(u_i)$ fix\'es par $F$. 
Remarquons que l'ordre de $A(u_i)$ ne d\'epend pas de $i$. 

\smallskip
Lusztig a conjectur\'e dans \cite{LMadison}, puis d\'emontr\'e dans
\cite{LUS}, qu'\'etant donn\'e un caract\`ere irr\'eductible $\rho$ de
$G(\mathbb{F}_q)$, il existe une unique classe unipotente
$\orb_\rho^{\Lu}$ de $G$ telle que $\VM(\rho,\orb_\rho^{\Lu})\ne 0$ et
qui soit de dimension maximale pour cette propri\'et\'e (sous
l'hypoth\`ese que la caract\'eristique $p$ de $\mathbb{F}_q$ est
suffisamment grande). En utilisant \cite{LUS}, Geck et Malle ont
montr\'e dans \cite{GMU} qu'il existe une unique classe unipotente
$\orb^{\GM}_\rho$ de $G$ telle que $\VMP(\rho,\orb^{\GM}_\rho)\ne 0$
et qui soit de dimension maximale pour cette propri\'et\'e (sans
restriction sur $p$). Les classes $\orb_\rho^{\Lu}$ et
$\orb_\rho^{\GM}$ \'etant \'egales, d'apr\`es \cite[Theorem 1.4]{GVM},
nous les noterons d\'esormais simplement $\orb_\rho$.

\smallskip
D'apr\`es \cite[Theorem~3.7]{GMU}, si le caract\`ere $\rho$ appartient \`a
$\E(G^F)_{s,\bc}$, la classe $\orb_\rho$ est \'egale \`a la classe
$\orb_{s,\bc}$. La d\'efinition de $\VM(\rho,\orb)$ montre d'autre part
que si la restriction de la valeur moyenne d'un caract\`ere \`a une
classe unipotente est non nulle, alors la restriction du caract\`ere
lui-m\^eme \`a cette classe est aussi non nulle. Par cons\'equent, la
restriction de $\rho$ \`a la classe $\orb_{s,\bc}$ est non nulle.

\begin{remark}
En mauvaise caract\'eristique, il peut arriver que 
le support unipotent de $\rho$ (au sens de l'introduction) soit
diff\'erent de $\orb_\rho$. Par exemple, dans le cas o\`u $G$ est un
groupe simple de type $G_2$, o\`u $p$ est \'egal \`a $3$ et o\`u $\orb$
est la classe des \'el\'ements unipotents r\'eguliers, il existe des
caract\`eres unipotents de $G^F$ dont la restriction \`a $\orb^F$ est non
identiquement nulle alors que leur valeur moyenne sur $\orb^F$ est nulle
(voir \cite{E}).
\end{remark} 

\subsection{Supports de valeurs moyennes de caract\`eres}

Dans cette sous-section, nous n'imposons pas de condition sur $p$ (\ie le
cas $p$ mauvais est permis).
Pour tout \'el\'ement $w$ de $W$, nous notons $T_w$ le tore maximal
rationnel obtenu \`a partir de $T$ par torsion par l'\'el\'ement $w$.
Pour tout $E_1\in\Irr(W_s)$, nous d\'efinissons la combinaison
lin\'eaire rationnelle suivante de caract\`eres de Deligne-Lusztig de
$G_s$: 
\begin{equation} \label{eqn:RsE}
R_s(E_1)=|W_s|^{-1}\,\sum_{w\in W_s}\Tr(w,E_1)\,R_{T_w}^{G_s}(1)
.\end{equation}
La fonction centrale $R_s(E_1)$ sur $G^F$ ainsi d\'efinie est celle
de \cite[(3.7.1)]{Lbook} et co\"\i ncide avec le caract\`ere fant\^ome
$R_{x_{E_1}}$.

Le r\'esultat suivant est prouv\'e dans \cite[Proposition 4.3]{GMU} pour
les caract\`eres unipotents. Nous l'\'etendons ici \`a tous les
caract\`eres irr\'eductibles.

\begin{thm}
Soit $\rho$ un caract\`ere irr\'eductible du groupe $G(\mathbb{F}_q)$ des
points $\mathbb{F}_q$-rationnels d'un groupe r\'eductif connexe $G$. 
Toute classe
unipotente rationnelle $\orb$ telle que $\VMP(\rho,\orb)\ne 0$ (resp.
$\VM(\rho,\orb)\ne 0$) est contenue dans l'adh\'erence de Zariski de $\orb_\rho$.
\end{thm}
\begin{proof}
Nous notons $V(\rho,\orb)$ soit $\VMP(\rho,\orb)$ soit 
$\VM(\rho,\orb)$.
Nous commen\c cons par nous ramener au cas o\`u le centre du groupe $G$ est
connexe.
Pour cela, nous fixons un plongement r\'egulier (voir \cite[chap. 14]{Lbook} ou
\cite{Ldisc}) de $G$ dans un groupe $G_0$ \`a centre connexe de m\^eme
groupe d\'eriv\'e que $G$ et un caract\`ere irr\'eductible $\rho_0$
de $G_0(\mathbb{F}_q)$ tel que $\rho$ intervienne dans la restriction de
$\rho_0$ \`a $G(\mathbb{F}_q)$. Notons $m_\rho$ ($\ge 1$) le nombre de
caract\`eres irr\'eductibles de $G(\mathbb{F}_q)$ qui apparaissent dans la 
restriction de $\rho_0$ \`a $G_0(\mathbb{F}_q)$ (ce nombre est
ind\'ependant du choix de $\rho_0$). On a alors (voir \cite[preuve du
Theorem~3.7]{GMU}) 
$$|A_{G_0}(\orb)|\,V(\orb,\rho_0)=m_\rho|A_G(\orb)|\,V(\orb,\rho).$$ 
Il suffit donc de prouver le th\'eor\`eme pour $G_0$. 
Nous supposons d\'esormais que le centre de $G$ est connexe.

Il existe un \'el\'ement semi-simple $s$ du dual de Langlands de $G$ et une
repr\'esentation $E_1$ du groupe de Weyl $W_s$ de $\Cent_{G^*}(s)$ tel que
le produit scalaire de $\rho$ et du caract\`ere fant\^ome $R_s(E_1)$
(d\'efini en~\ref{eqn:RsE}) soit non nul. 
D'autre part, toutes les repr\'esentations irr\'eductibles $E'$ de $W_s$, dont
les caract\`eres fant\^omes $R_s(E')$ associ\'es ont produit scalaire non nul
avec $\rho$ appartiennent \`a une m\^eme cellule bilat\`ere $\bc_1$ de $W_s$.
La projection uniforme $\rho_{\unif}$ de $\rho$ ({\it i.e.,} la projection
sur l'espace des combinaisons lin\'eaires de caract\`eres de Deligne-Lusztig) 
est combinaison lin\'eaire de caract\`eres fant\^omes $R_s(E')$ tels 
que $E'$ appartienne \`a $\bc_1$:
$$\rho_{\unif}=\sum_{E'\in\bc_1}c(\rho,E')R_s(E').$$
La valeur moyenne pond\'er\'ee $V(\orb,\rho)$ \'etant, d'apr\`es 
\cite[Proposition 1.3]{GVM}, \'egale au produit scalaire de $\rho$ par une 
fonction uniforme ({\it i.e.,} combinaison lin\'eaire de caract\`eres de 
Deligne-Lusztig), on a $V(\orb,\rho)=V(\orb,\rho_{\unif})$. 
Par cons\'equent,
$$V(\orb,\rho)=\sum_{E'\in\bc_1}c(\rho,E')\,V(\orb,R_s(E')).$$
Soit $\bc$ la cellule bilat\`ere de $W$ induite de $\bc_1$.
En utilisant la Proposition~\ref{prop:ind-supp}, nous voyons que 
$V(\orb,\rho)$ est combinaison lin\'eaire de valeurs moyennes
pond\'er\'ees $V(\orb,R_s(E))$, o\`u $E$ appartient \`a une cellule
bilat\`ere $\bc'$ telle que $\bc'\le \bc$ et il suffit alors
de raisonner comme dans la d\'emonstration de \cite[Proposition 4.3]{GMU}. 
\end{proof}

\end{document}